\documentclass{article}
\usepackage{amsmath,amsxtra,amssymb,latexsym,amsfonts,amscd}
\usepackage{multicol,color}
\usepackage{float}
\usepackage{soul}
\usepackage{graphicx}
\usepackage{pstricks-add}
\usepackage{pgf,tikz}
\usepackage{mathrsfs}
\usetikzlibrary{arrows}
\usepackage{amssymb}
\usepackage{theorem}
\usepackage{fancyhdr}
\usepackage{tikz}
\usepackage{enumerate}
\usepackage[margin=1.1in]{geometry}
\def\disp{\displaystyle}
\def\e{\epsilon}

\def\dd{\delta}
\def\DD{\Delta}
\def\lm{\lambda}
\def\O{\Omega}
\def\Tilde{\widetilde}
\def\tilde{\widetilde}

\def\({\left(}
\def\){\right)}
\def\[{\left[}
\def\]{\right]}
\def\n{\left \|}
\def\en{\right \|}
\def\nn{\left \{ }
\def\hnn{\right \}}

\def\ox{\bar{x}}
\def\oy{\bar{y}}
\def\oz{\bar{z}}
\def\ov{\bar{v}}
\def\ou{\bar{u}}
\def\ot{\bar{t}}

\def\gph{\hbox{}}

\def\gg{\gamma}
\def\dn{\downarrow}

\def\tto{\rightrightarrows}

\def\tilde{\widetilde}
\def\Tilde{\widetilde}
\def\Bar{\overline}
\def\la{\langle}
\def\ra{\rangle}

\def\h{\hfill\Box}
\def\R{\mathbb{R}}
\def\N{\mathbb{N}}

\def\co{\mbox{\rm co}\,}
\def\gph{\mbox{\rm gph}\,}
\def\epi{\mbox{\rm epi}\,}

\def\dom{\mbox{\rm dom}\,}

\def\proj{\mbox{\rm proj}\,}

\def\dist{\mbox{\rm dist}}

\def\dn{\downarrow}
\def\O{\Omega}
\def\o{\omega}

\def\vph{\varphi}
\def\emp{\emptyset}

\def\oR{\Bar{\R}}
\def\lm{\lambda}

\def\gg{\gamma}

\def\dd{\delta}
\def\DD{\Delta}
\def\al{\alpha}
\def\vth{\vartheta}

\def\be{\beta}
\def\vth{\vartheta}
\def\ph{\varphi}

\def\N{I\!\!N}
\def\th{\theta}

\newtheorem{theorem}{Theorem}[section]
\newtheorem{lemma}[theorem]{Lemma}
\newtheorem{corollary}[theorem]{Corollary}
\newtheorem{proposition}[theorem]{Proposition}
\newtheorem{definition}[theorem]{Definition}
\theoremstyle{plain}{\theorembodyfont{\rmfamily}
}
\theoremstyle{plain}{\theorembodyfont{\rmfamily}
}
\theoremstyle{plain}{\theorembodyfont{\rmfamily}
}
\theoremstyle{plain}{\theorembodyfont{\rmfamily}
\newtheorem{example}[theorem]{Example}}

\theoremstyle{plain}{\theorembodyfont{\rmfamily}
}

\def\eq{\begin{equation}}
\def\eeq{\end{equation}}
\begin{document}
\begin{center}
{\bf OPTIMIZATION OF A PERTURBED SWEEPING PROCESS\\BY CONSTRAINED DISCONTINUOUS CONTROLS}\\[3ex]
GIOVANNI COLOMBO\footnote{Dipartimento di Matematica ``Tullio Levi-Civita", Universit$\grave{\textrm{a}}$ di Padova,
via Trieste 63, 35121 Padova, Italy (colombo@math.unipd.it) and G.N.A.M.P.A. of INdAM.} \quad BORIS S. MORDUKHOVICH\footnote{Department of Mathematics, Wayne State University, Detroit, Michigan 48202, USA (boris@math.wayne.edu). Research of this author was partly supported by the USA National Science Foundation under grants DMS-1512846 and DMS-1808978, and by the USA Air Force Office of Scientific Research grant \#15RT0462.} \quad DAO NGUYEN\footnote{Department of Mathematics, Wayne State University, Detroit, Michigan 48202, USA (dao.nguyen2@wayne.edu). Research of this author was partly supported by the USA National Science Foundation under grant DMS-1808978 and by the USA Air Force Office of
Scientific Research grant \#15RT0462.}
\end{center}
\small{\sc Abstract.} This paper deals with optimal control problems described by a controlled version of Moreau's sweeping process governed by convex polyhedra, where measurable control actions enter additive perturbations. This class of problems, which addresses unbounded discontinuous differential inclusions with intrinsic state constraints, is truly challenging and underinvestigated in control theory while being highly important for various applications. To attack such problems with constrained measurable controls, we develop a refined method of discrete approximations with establishing its well-posedness and strong convergence. This approach, married to advanced tools of first-order and second-order variational analysis and generalized differentiations, allows us to derive adequate collections of necessary optimality conditions for local minimizers, first in discrete-time problems and then in the original continuous-time controlled sweeping process by passing to the limit. The new results include an appropriate maximum condition and significantly extend the previous ones obtained under essentially more restrictive assumptions. We compare them with other versions of the maximum principle for controlled sweeping processes that have been recently established for global minimizers in problems with smooth sweeping sets by using different techniques. The obtained necessary optimality conditions are
illustrated by several examples.\\[1ex]
\vspace*{-0.2in}

\normalsize
\section{Introduction and Problem Formulation}\label{assumptions}

This paper addresses the following optimal control problem labeled as $(P)$:\\[1ex]
Minimize the Mayer-type cost functional
\begin{equation}\label{cost1}
J[x,u]:=\vph\big(x(T)\big)
\end{equation}
over the corresponding (described below) pairs $(x(\cdot),u(\cdot))$ satisfying
\begin{equation}\label{Problem}
\left\{\begin{matrix}
\dot{x}(t)\in-N\big(x(t);C\big)+g\big(x(t),u(t)\big)\;\textrm{ a.e. }\;t\in[0,T],\;x(0)=x_0\in C\subset\R^n,\\
u(t)\in U\subset\R^d\;\textrm{ a.e. }\;t\in[0,T],
\end{matrix}\right.
\end{equation}
where the set $C$ is a convex {\em polyhedron} given by
\begin{equation}\label{C}
C:=\bigcap_{j=1}^{s}C^j\textrm{ with }C^j:=\nn x\in\R^n\big|\;\la x^j_*,x\ra\le c_j\hnn,
\end{equation}
and where $N(x;C)$ stands for the normal cone of convex analysis defined by
\begin{equation}\label{nc}
N(x;C):=\big\{v\in\R^n\;\big|\;\la v,y-x\ra\le 0,\;y\in C\big\}\textrm{ if }x\in C\textrm{ and }N(x;C):=\emp\textrm{ if }x\notin C.
\end{equation}
Observe that due to the second part of definition \eqref{nc} mandatory yields the presence of the hidden
{\em pointwise state constraints} on the trajectories of \eqref{Problem}:
\begin{equation}\label{e:8}
x(t)\in C,\textrm{ i.e. }\la x^j_*,x(t)\ra\le c_j\;\textrm{ for all }\;t\in [0,T]\;\textrm{ and }\;j=1,\ldots,s.
\end{equation}
Considering the differential inclusion in \eqref{Problem} without the additive perturbation term $g(x,u)$, we arrive at
the framework of the {\em sweeping process} introduced by Jean-Jacques Moreau who was motivated by applications
to problems of elastoplasticity; see \cite{mor_frict}. It has been well recognized that the (uncontrolled) Moreau's sweeping
process has a {\em unique} absolutely continuous (or even Lipschitz continuous) solution for convex and mildly nonconvex
sets $C$; see, e.g., \cite{CT} and the references therein. Thus there is no room for optimization of the sweeping process
unless some additional functions or parameters of choice are inserted into its description. It is very different from
control theory for Lipschitzian differential inclusions
\begin{equation}\label{diff-inc}
\dot x(t)\in F\big(x(t)\big)\textrm{ a.e. }t\in[0,T],\;x(0)=x_0\in\R^n,
\end{equation}
which have multiple solutions. The latter type of dynamics extends the classical ODE control setting with $F(x):=f(x,U)$ in \eqref{diff-inc}, where
the choice of measurable controls $u(t)\in U\subset\R^d$ a.e.\ $t\in[0,T]$ creates the possibility to find an optimal one with respect to a
prescribed performance. The main issue here is that the normal cone mapping $N(\cdot;C)$ in the sweeping process is {\em highly
non-Lipschitzian} (even discontinuous) while being {\em maximal monotone}. On the other hand, the well-developed optimal control theory for
differential inclusions \eqref{diff-inc} strongly depends on Lipschitzian behavior of $F(\cdot)$; see, e.g., \cite{m-book2,v} with the references
therein as well as more recent publications.

Introducing controls into the perturbation term of \eqref{Problem} allows us to have multiple solutions $x(\cdot)$ of this system by the choice of
feasible control functions $u(\cdot)$ and thus to minimize the cost functional \eqref{cost1} over feasible control-trajectory pairs. Problems of this
type were considered in the literature from the viewpoint of the existence of optimal solutions and relaxation; see \cite{aht,cmf,et,Tol} among
other publications. 

More recently, {\em necessary optimality conditions} for local minimizers were derived in \cite{cm1,cm2} by the {\em method of discrete
approximations} for problems of type $(P)$ with smooth (in fact $W^{2,\infty}$) control functions without any constraints. Later on these results
were further extended in \cite{cm3} to nonconvex (and hence nonpolyhedral) problems with prox-regular sets $C$ in the same control setting.
Note that both $C$ and $g$ in \eqref{Problem} may be time-dependent; we discuss the autonomous case just for simplicity. The discrete
approximation approach implemented in \cite{cm1}--\cite{cm3} was based on the scheme from \cite{chhm1} developed for the unperturbed
sweeping process with controls in the moving set. The later was in turn a sweeping control version of the original discrete approximations
method to derive necessary optimality conditions for Lipschitzian differential inclusions \eqref{diff-inc} suggested and implemented in \cite{m95};
see also \cite{m-book2}.

Quite recently, other approximation procedures were developed to derive necessary optimality conditions for global minimizers of $(P)$ in the
class of measurable controls while under rather strong assumptions. The first paper \cite{ac} assumes, among other requirements, that the
boundary of the sweeping set $C$ in \eqref{Problem} is ${\cal C}^3$-smooth, the control set $U$ is compact and convex, and its image $g(x,U)$
under $g$ is convex as well. The ${\cal C}^3$-smoothness assumption on $C$ was relaxed in \cite{pfs}, by employing a smooth approximation
procedure not relying on the distance function as in \cite{ac}, for the case of $C:=\{x\in\R^n\;|\;\psi(x)\le 0\}$ with $\psi$ being a ${\cal C}^2$-
smooth convex function. The necessary optimality conditions obtained in both papers \cite{ac,pfs} can be treated as somewhat different
counterparts of the celebrated Pontryagin Maximum Principle (PMP) for state-constrained controlled differential equations $\dot x=f(x,u)$.

Note that necessary optimality conditions in some other classes of optimal control problems governed by various controlled versions of the
sweeping process were developed in \cite{ao,bk,cm1,cm2,cm3,chhm,chhm1,hm18}.

The main goal of this paper is to derive necessary optimality conditions for local minimizers (in the senses specified below) of the formulated
problem $(P)$, with the constraint set $U$ in \eqref{Problem} given by an {\em arbitrary compact} and with the (nonsmooth) {\em polyhedral} set $C$ from \eqref{C}, by significantly reducing regularity assumptions on the reference control. Although problem \eqref{Problem} is stated in the class of {\em measurable} feasible control actions, we assume that the local {\em optimal} control under consideration is of {\em bounded variation}, hence allowing to be discontinuous.

Our approach is based on developing the {\em method of discrete approximations}, which is certainly of its own interest and has never been implemented before in control theory for sweeping processes with discontinuous controls. The novel results in this direction establish a {\em strong} approximation of {\em every feasible} control-state pair for $(P)$ in the sense of the $L^2$-norm convergence of discretized controls and the $W^{1,2}$-norm convergence of the corresponding piecewise linear trajectories.
Furthermore, we justify such a strong convergence of {\em optimal} solutions for discrete problems to the given local minimizer of $(P)$.

Dealing further with {\em intrinsically nonsmooth} and {\em nonconvex} discrete-time approximation problems, we derive for
them necessary optimality conditions of the discrete Euler-Lagrange type by using appropriate unconvexified tools of first-order and
second-order variational analysis and generalized differentiation. Employing these tools and passing to the limit from discrete approximations lead us to new nondegenerate necessary optimality conditions for local optimal solutions of the sweeping control problem $(P)$. The obtained results
significantly extend those recently established in \cite{cm2} for unconstrained $W^{2,\infty}$ optimal controls in $(P)$, contain a maximum
condition, while being essentially different from the necessary optimality conditions derived in \cite{ac,pfs} for problems of type $(P)$ with
smooth sets $C$ in addition to other assumptions. We present nontrivial examples that illustrate the efficiency of the new results. Further
applications to some practical models are considered in our subsequent paper \cite{cmn}.

The rest of the paper is organized as follows. In Section~2 we formulate the standing assumptions, discuss the types of local minimizers under
consideration, and present some preliminary results.

Section~3 is devoted to the construction of discrete approximations of the controlled constrained sweeping dynamics \eqref{Problem} that
allows us to deal with measurable controls (in fact of bounded variation) and to strongly approximate any feasible solutions of $(P)$ as mentioned above. This result plays a major role in the justification of the developed version the method of discrete approximations for problem $(P)$.

In Section~4 we construct a sequence of discrete approximation of a given ``intermediate" local minimizer for $(P)$ that occupies an
intermediate position between weak and strong minimizers in variational and control problems. The major result of this section justifies the
strong $ W^{1,2}\times L^2$ approximation of the given local minimum pair $(\ox(\cdot),\ou(\cdot))$ by extended optimal solutions to the
discretized problems. It makes a bridge between the continuous-time sweeping control problem $(P)$ and its discrete-time counterparts.

It occurs that the discrete-time approximating problems are unavoidably nonsmooth and nonconvex, even when the initial data are
differentiable. It is due to the presence of increasingly  many geometric constraints generated by the normal cone graph. To deal with them, we
need adequate tools of variational analysis involving not only first-order but also second-order generalized differentiation. The latter is because
of the normal cone description of the sweeping process. In Section~5 we present the corresponding definitions of the first-order and second-
order generalized differential constructions taken from \cite{m-book1} together with the results of their computations entirely in terms of the given
data of \eqref{Problem}.

Section~6 provides the derivation of necessary optimality conditions for discrete-time problems by reducing them to problems of
nondifferentiable programming with many geometric constraints, using necessary optimality conditions for them obtained via variational/
extremal principles, and then expressing the latter in terms of the given data of $(P)$ by employing calculus rules of generalized differentiation.

Section~7 is the culmination. We pass to the limit from the necessary optimality conditions for discrete-time problems by using stability of
discrete approximations, robustness of then generalized differential constructions, and establishing an appropriate convergence of adjoint
functions, which is the most difficult part. In this way we arrive at new necessary conditions for local minimizers of $(P)$ expressed in terms of
the given data of the original problem. The usefulness of the nondegenerated optimality conditions obtained is illustrated in Section~8 by nontrivial examples.

Throughout the paper we use standard notations of variational analysis and optimal control; see, e.g., \cite{m-book1,m-book2}. Recall that
$\mathbb{N}:=\{1,2,\ldots\}$.\vspace*{-0.2in}

\section{Standing Assumptions and Basic Notions}\label{stand}
\setcounter{equation}{0}\vspace*{-0.1in}

Dealing with the polyhedron $C$ from \eqref{C} and having $\ox\in C$, consider the set of {\em active constraint indices}
\begin{equation}\label{aci}
I(\ox):=\big\{j\in\{1,\ldots,s\}\;\big|\;\la x^j_*,\ox\ra=c_j\big\}.
\end{equation}
Recall that the {\em linear independence constraint qualification} (LICQ) holds at $\ox$ if
\begin{equation}\label{plicq}
\Big[\sum_{j\in I(\ox)}\al_jx^j_*=0,\;\al_j\in\R\Big]\Longrightarrow\big[\al_j=0\textrm{ for all }j\in I(\ox)\big\}.
\end{equation}

Our {\em standing assumptions} in this paper are as follows:\\[1ex]
{\bf(H1)} The control region $U\ne\emp$ is a compact set in $\R^d$ (in fact it may be an arbitrary metric compact).\\
{\bf(H2)} The perturbation mapping $g\colon\R^n\times U\to\R^n$ is continuous in $(x,u)$ while being also Lipschitz continuous with respect to
$x$ uniformly on $U$ whenever $x$ belongs to a bounded subset of $\R^n$ and satisfies there the sublinear growth condition
\begin{equation*}
\|g(x,u)\|\le\be\big(1+\|x\|\big)\;\mbox{ for all }\;u\in U
\end{equation*}
with some positive constant $\be$.\\
{\bf(H3)} The LICQ condition \eqref{plicq} holds along the reference trajectory $\ox(t)$ of \eqref{Problem} for all $t\in[0,T]$.

It follows from \cite[Theorem~1]{et} that for each measurable control $u(\cdot)$ there is a unique solution $x(\cdot)\in W^{1,2}([0,T],\R^n)$ to the
Cauchy problem in \eqref{Problem}. Thus by a {\em feasible process} for $(P)$ we understand a pair $(x(\cdot),u(\cdot))$ such that $u(\cdot)$ is
measurable, $x(\cdot)\in W^{1,2}([0,T],\R^n)$, and all the constraints in \eqref{Problem} are satisfied. The above discussion tells us that the set
of feasible pairs for $(P)$ is nonempty.

Furthermore, it follows from \cite[Theorem~2]{et} that under the assumptions above the sweeping control problem $(P)$ admits an {\em optimal
solution} provided that the image set
\begin{equation*}
g(x,U):=\big\{x\in\R^n\;\big|\;x=g(x,u)\textrm{ for some }u\in U\big\}
\end{equation*}
is convex. Since in this paper we are interested in deriving necessary optimality conditions for a given local minimizer of $(P)$, we do not
impose the aforementioned convexity assumption.

Let us now specify what we mean by a local minimizer of $(P)$.\vspace*{-0.07in}

\begin{definition}\label{Def3.1}
We say that a feasible pair $(\ox(\cdot),\ou(\cdot))$ for $(P)$ is a {\sc $W^{1,2}\times L^2$-local minimizer} in this problem if $\ox(\cdot)\in
W^{1,2}([0,T];\R^n)$ and there exists $\e>0$ such that $J[\ox,\ou]\le J[x,u]$ for all feasible pairs $(x(\cdot),u(\cdot))$ satisfying the condition
\begin{equation*}
\int_0^T\(\n\dot{x}(t)-\dot{\ox}(t)\en^2+\n u(t)-\ou(t)\en^2\)dt<\e.
\end{equation*}
\end{definition}\vspace*{-0.05in}

For the case of differential inclusions of type \eqref{diff-inc} with no explicit controls, this notion corresponds to {\em intermediate local
minimizers} of rank two introduced in \cite{m95} and then studied there and in other publications; see, e.g., \cite{m-book2,v} and the references
therein. Quite recently, such minimizers have been revisited in \cite{hm18} for controlled sweeping processes different from \eqref{Problem};
namely, for those where continuous control actions enter the moving set $C(t)=C(u(t))$. It is easy to see that {\em strong} ${\cal C}\times L^2$-
local minimizers of $(P)$ with $\ox(\cdot)\in W^{1,2}([0,T];\R^n)$ fall into the category of Definition~\ref{Def3.1}, but not vice versa.

In the general setting of $W^{1,2}\times L^2$-local minimizers we need to use a certain relaxation procedure in the line of Bogolyubov and
Young that has been well understood in the calculus of variations and optimal control; see, e.g., \cite{dfm,et,m-book2,Tol,v} for more recent
publications in the case of differential inclusions. Taking into account the convexity and closedness of the normal cone $N(x;C)$ and the
compactness of the set $g(x,U)$, the {\em relaxed} version $(R)$ of problem $(P)$ consists of minimizing the cost functional \eqref{cost1} on
absolutely continuous trajectories of the convexified differential inclusion
\begin{equation}\label{conv}
\dot{x}(t)\in-N\big(x(t);C\big)+\co g\big(x(t),U\big)\;\textrm{ a.e. }\;t\in[0,T],\;x(0)=x_0\in C\subset\R^n,
\end{equation}
where `co' signifies the convex hull of the set. Then we come up with the following notion.\vspace*{-0.07in}

\begin{definition}\label{relaxed} Let $(\ox(\cdot),\ou(\cdot))$ be a feasible pair for $(P)$. We say that it is a {\sc relaxed $W^{1,2}\times
L^2$-local minimizer} for $(P)$ if $\ox(\cdot)\in W^{1,2}([0,T];\R^n)$ and there is $\e>0$ such that
\begin{equation*}
\ph\big(\ox(T)\big)\le\ph\big(x(T)\big)\;\textrm{ whenever }\;\int_0^T\(\n\dot x(t)-\dot{\ox}(t)\en^2+\n u(t)-\ou(t)\en^2\)dt<\e,
\end{equation*}
where $u(\cdot)$ is a measurable control with $u(t)\in\co U(t)$ a.e.\ on $[0,T]$, and where $x(\cdot)$ is a trajectory of the convexified inclusion
\eqref{conv}
that can be strongly approximated in $W^{1,2}([0,T];\R^n)$ by feasible trajectories to $(P)$ generated by piecewise constant controls
$u_m (\cdot)$ on $[0,T]$ with
\begin{equation*}
\int_0^T\|u_m(t)-u(t)\|^2dt\to 0\;\textrm{ as }\;m\to\infty.
\end{equation*}
\end{definition}\vspace*{-0.05in}

Since step functions are dense in the space $L^2([0,T];\R^d)$, we obviously have that there is no difference between $W^{1,2}\times L^2$-local
minimizers for $(P)$ and their relaxed counterparts provided that the sets $g(x,U)$ and $U$ are {\em convex}, which is not assumed in what
follows. Moreover, it is possible to deduce from the proofs of \cite[Theorem~2]{et} and \cite[Theorem~4.2]{Tol} that any strong local minimizer for
$(P)$ is automatically a relaxed one under the assumptions made, but we are not going to pursue this issue here.

Consider further a set-valued mapping $F\colon C\times U\tto\R^n$ defined by
\begin{equation}\label{F0}
F(x,u):=N(x;C)-g(x,u)\textrm{ for all }x\in C,\;u\in U
\end{equation}
and deduce from the Motzkin's theorem of the alternative the representation
\begin{equation}\label{F}
F(x,u):= \Big\{\sum_{j\in I(x)}\lm^j x^j_*\;\Big|\;\lm^j\ge 0\Big\}-g(x,u),\quad x\in C,\;u\in U,
\end{equation}\vspace*{-0.35in}

\section{Discrete Approximations of Feasible Solutions}\label{disc-giov}
\setcounter{equation}{0}\vspace*{-0.1in}

In this section we start developing the {\em method of discrete approximations} to study the sweeping control problem $(P)$ under our standing
assumptions. For simplicity, consider the standard Euler explicit scheme for the replacement of the time derivative in \eqref{Problem} by
\begin{equation*}
\dot x(t)\approx\frac{x(t+h)-x(t)}{h}\textrm{ as }h\dn 0,
\end{equation*}
which we formalize as follows. For any $m\in\N$ denote by
\begin{equation*}
\DD_m:=\big\{0=t^0_m<t^1_m<\ldots<t^{2^m}_m=T\big\}\textrm{ with }h_m:=t^{i+1}_m-t^i_m
\end{equation*}
the discrete mesh on $[0,T]$ and define the sequence of discrete-time systems
\begin{equation}\label{e:3.4}
x^{i+1}_m\in x^i_m-h_m F(x_m^{i},u_m^{i}),\;i=0,\ldots,2^m-1,
\end{equation}
where we have $u_m^{i}\in U$ due to the definition of $F$ in \eqref{F0}. Let $I_m^i:=[t_m^{i-1},t_m^i)$.

The next result provides a constructive approximation of {\em any} feasible process for $(P)$ by feasible solutions to \eqref{e:3.4} that are
appropriately extended to the continuous-time interval $[0,T]$. This result plays a major role in the entire subsequent procedure to derive
necessary optimality conditions for $(P)$ while certainly being of its independent interest. Recall that a {\em representative} of a given
measurable function on $[0,T]$ is a function that agrees with the given one for a.e.\ $t\in[0,T]$.\vspace*{-0.07in}

\begin{theorem}\label{Thm3.1} Let $(\ox(\cdot),\ou(\cdot))$ be a feasible pair for problem $(P)$ such that $\ox(\cdot)\in W^{1,2}([0,T];\R^n)$ and
that $\ou(\cdot)$ is of bounded variation $($BV$)$ while admitting a right continuous representative on $[0,T]$, which we keep denoting by
$\ou(\cdot)$. In addition to {\rm(H1)--(H3)}, suppose that the mapping $g(x,u)$ is locally Lipschitzian in both variables around
$(\ox(t),\ou(t))$ for all  $t\in[0,T]$. Then for each $i=1,\ldots,2^m$ there exist sequences of unit vectors $z^{ji}_m\in\R^n$,
real numbers $c^{ji}_m$, and state-control  pairs $(x_m(t),u_m(t))$, $0\le t\le T$, such that
\begin{equation}\label{3.1a}
z^{ji}_m\to x^j_*\textrm{ and }c^{ji}_m\to c_j\textrm{ as }\;m\to\infty
\end{equation}
and the following properties are fulfilled:\\
{\bf (a)} The sequence of control mappings $u_m\colon[0,T]\to U$, which are constant on each interval $I_m^i$, converges to $\ou(\cdot)$
strongly in  $L^2([0,T];\R^d)$ and pointwise on $[0,T]$.\\
{\bf (b)} The sequence of continuous state mappings $x_m\colon[0,T]\to\R^n$, which are affine on each interval $I_m^i$, converges strongly in
$W^{1,2}([0,T];\R^n)$ to $\ox(\cdot)$ while satisfying the inclusions
\begin{equation}\label{3.1b}
x_m(t_m^i)=\ox(t_m^i)\in C_m^i\textrm{ for each }i=1,\ldots,2^m\textrm{ with }x_m(0)=x_0,
\end{equation}
where the perturbed polyhedra $C_m^i$ are given by
\begin{equation}\label{3.1C}
C_m^i:=\bigcap_{j=1}^s\big\{x\in\R^n\;\big|\;\langle z_m^{ji},x\rangle\le c_m^{ji}\big\}\textrm{ for }i=1,\ldots,2^m\textrm{ with }C_m^0:=C.
\end{equation}
{\bf (c)} For all $t\in(t_m^{i-1},t_m^i)$ and $i=1,\ldots,2^m$ we have the differential inclusions
\begin{equation}\label{3.1c}
\dot{x}_m(t)\in-N\big(x_m(t_m^{i});C_m^i\big)+g\big(x_m(t_m^{i}),u_m(t)\big).
\end{equation}
\end{theorem}\vspace*{-0.05in}

As a part of the proof of Theorem~\ref{Thm3.1}, we establish the following lemma, which is of its own interest.\vspace*{-0.07in}

\begin{lemma}\label{lemma1} Given a feasible solution $(\ox(\cdot),\ou(\cdot))$ to $(P)$ under the assumptions of Theorem~{\rm\ref{Thm3.1}},
we have:

{\bf (i)} $\ox(\cdot)$ is Lipschitz continuous on $[0,T]$ and right differentiable for every $t\in[0,T]$, and its right derivative
denoted by $\dot\ox(\cdot)$ is also right continuous on $[0,T]$.

{\bf (ii)} The sweeping differential inclusion
\begin{equation*}
\dot{\ox}(t)\in-N\big(\ox(t);C\big)+g\big(\ox(t),\ou(t)\big),
\end{equation*}
with $\dot{\ox}(t)$ taken from {\rm(i)} and the right continuous representative of $\ou(t)$, is satisfied for each $t\in[0,T]$.
\end{lemma}\vspace*{-0.07in}
{\bf Proof.} Considering the differential inclusion
\begin{equation*}
\dot{x}(t)\in-N(x(t);C)+g\big(x(t),\ou(t)\big),\;x(0)=x_0\in C,
\end{equation*}
we deduce from, e.g., \cite[Propositions~3.8 and 3.12]{Bre} that there exists one and only one Lipschitz continuous solutions on $[0,T]$, which
therefore agrees with the given trajectory $\ox(\cdot)$. Furthermore, $\ox(\cdot)$ is a unique solution of the differential inclusion
\begin{equation*}
\dot{x}(t)\in-N\big(x(t);C\big)+g\big(\ox(t),\ou(t)\big)\textrm{ a.e. }t\in[0,T],\;x(0)=x_0\in C,
\end{equation*}
where the perturbation term depends only on $t$. Then the assumptions imposed on $\ou(\cdot)$ and $g$ ensure that the mapping
$t\mapsto g(\ox(t),\ou(t))$ is BV on $[0,T]$. The result of \cite[Proposition~3.3]{Bre} tells us that $\ox(\cdot)$ is right differentiable at each
$t\in[0,T)$ and satisfies the equalities
\begin{equation}\label{i2}
\dot{\ox}(t)=g\big(\ox(t),\ou(t)\big)-\proj_{N(\ox(t);C)}\big(g(\ox(t),\ou(t)\big)=\proj_{T(\ox(t);C)}\big(g(\ox(t),\ou(t)\big))\textrm{ a.e. }t\in[0,T]
\end{equation}
written via the (unique) projection onto the convex set $N(\ox(t);C)$, where the second one can be easily verified. Our goal is to show that $\dot{\ox}(t)$ is right continuous on $[0,T]$ while satisfying \eqref{i2} {\em for each} $t\in[0,T]$.

Observe preliminary that, thanks to LICQ, the polyhedron $C$ has nonempty interior and so the normal cone $N(x;C)$
is pointed at each $x\in\partial C$. Denote $v(t):=g(\ox(t),\ou(t))$, fix $0\le\ot<T$, and let $t_k\to\ot^+$ as $k\to\infty$.
We need to verify that $v(t_k)\to v(\ot)$, which is equivalent by \eqref{i2} to
\begin{equation}\label{i3}
\proj_{N(\ox(\ot);C)}\big(v(\ot)\big)=\underset{k\to\infty}{\lim}\proj_{N(\ox(t_k);C)}\big(v(t_k)\big).
\end{equation}
Note that there is nothing to prove if $\ox(\ot)\in{\rm int}\,C$, since $N(\ox(t_k);C)=\{0\}$ for all $k$ sufficiently large.
To proceed further, assume that $\ox(\ot)\in{\rm bd}\,C$ and observe easily that $N(\ox(t_k);C)\subset N(\ox(\ot);C)$
for all large $k$. Consider now the following three possible cases:

{\bf(1)} If $v(\ot)\in T(\ox(\ot);C)$, then $\langle v(\ot),x_{\ast}^j\rangle\le 0$ for all $j\in I(\ox(\ot))$, and so
$\lim_{k\to\infty}\langle v(t_k),x_{\ast}^j\rangle\le 0$ for all $j\in I(\ox(\ot))$. Since $N(\ox(t_k);C)\subset N(\ox(\ot));C)$, we get \eqref{i3}.

{\bf(2)} If $v(\ot)\in N(\ox(\ot);C)$, then arguing similarly to (1) and using the second equality in \eqref{i2} show that
$\mathrm{proj}_{T(\ox(t_k);C)}\big(v(t_k)\big)\to 0=\mathrm{proj}_{T(\ox(\ot);C)}(v(\ot))$ as $k\to\infty$ and hence verifies \eqref{i3} directly.

{\bf(3)} Let now $v(\ot)\not\in T(\ox(\ot);C)\cup N(\ox(\ot);C)$. Then LICQ ensures the unique representation
\begin{equation}\label{stella}
\mathrm{proj}_{N(\ox(\ot);C)}\big(v(\ot)\big)=\sum_{j\in I(\ox(\ot))}\lambda_j(\ot)x_{\ast}^j,
\end{equation}
where $\lambda_j(\ot)\ge 0$ for all $j\in I(\ox(\ot))$. If in this case $\lambda_j(\ot)=0$ for some $j\in I(\ox(\ot))$, then
$$
\big\langle\mathrm{proj}_{T(\ox(\ot);C)}\big(v(\ot)\big),x_{\ast}^j\big\rangle<0,
$$
due to the LICQ assumption. This implies that $\langle\ox(t_k),x_{\ast}^j\rangle\le\langle\ox(\ot),x_{\ast}^j\rangle$ for all $k$ sufficiently large.
Consequently, the corresponding vector $x_{\ast}^j$ appears being multiplied by zero in the representation
\begin{equation}\label{stella2}
\mathrm{proj}_{N(\ox(t_k);C)}\big(v(t_k)\big)=\sum_{j\in I(\ox(t_k))}\lambda_j(t_k)x_{\ast}^j.
\end{equation}
Recalling that $I(\ox(t_k))\subset I(\ox(\ot))$ for all large $k$, it turns out that the set of active indices in \eqref{stella2} is the same as in \eqref{stella}. Finally, it follows from \cite[Theorems~2.1 and 4.1]{rob} under the imposed LICQ that the coefficients $\lambda_j(\cdot)$ are continuous with respect to $v(\cdot)$. This verifies \eqref{i3} in case (3). $\h$\vspace*{0.02in}

Now we are ready to proceed with the proof of the major Theorem~\ref{Thm3.1}.\\[1ex]
{\bf Proof of Theorem~\ref{Thm3.1}}. Fix $m\in\mathbb{N}$ and for all $t\in[t_m^i,t_m^{i+1})$ and $i=0,\ldots,2^m-1$ define
\begin{equation*}
u_m(t):=\ou (t_m^{i+1}),\quad x_m(t):=\ox(t^i_m)+(t-t_m^i)\frac{\ox(t_m^{i+1})-\ox(t_m^i)}{h_m}.
\end{equation*}
Then denote $\o_m(t):=\dot\ox_m(t)$ for which we have the representation
\begin{equation*}
\o_m(t)=\o^i_m:=\frac{\ox(t^{i+1}_m)-\ox(t^i_m)}{h_m}\;\textrm{ whenever }\;t\in[t^i_m,t^{i+1}_m),\;i=0,\ldots,2^m-1.
\end{equation*}
It follows from the right continuity of $\ou$ that $u_m(t)\to\ou(t)$ as $m\to\infty$ for all $t\in[0,T)$. Hence we get that $u_m(\cdot)\to\ou(\cdot)$
strongly in $L^2(0,T)$ by the dominated convergence theorem, which verifies (a). To prove (b) and (c), let $\bar{t}$ be a nodal point of the $m$-th mesh that by construction remains a nodal point for all $m'$-mesh with $m'\ge m$. Denote by $i_m(\bar{t})$ the index $i$ such that $\bar{t}=i\frac{T}{2^m}$ and observe by Lemma~\ref{lemma1} that
\begin{equation}\label{lim3.1}
\lim_{m\to\infty}\o^{i_m(\bar{t})}_m=\dot{\ox}(\bar{t})\;\textrm{ and }\;\lim_{m\to\infty}\|\o_m-\dot{\ox}\|_{L^2(0,T)}=0.
\end{equation}
Indeed, the first equality in \eqref{lim3.1} is a consequence of the construction of $\o_m(t)$ and the right continuity of the derivative $\dot\ox(t)$ on $[0,T]$ by Lemma~\ref{lemma1}(i). This in turn yields the second equality therein by basic real analysis and thus justifies the strong $W^{1,2}$ convergence of $x_m(\cdot)$ to $\ox(\cdot)$.

To proceed further, for any $t\in[0,T)$ denote $\tau_m(t):=\min\{t_m^i\;|\;t\le t_m^i\}$ and get by assertion (a) that
\begin{eqnarray*}
\dot{x}_m(t)-g(x_m\big(\tau_m(t)),u_m(t)\big)=\o_m\big(\tau_m(t)\big)-g\big(x_m(\tau_m(t)),u_m(t)\big)
\rightarrow\dot{\ox}(t)-g\big(\ox(t),\ou(t)\big)\in-N\big(\ox(t);C\big)
\end{eqnarray*}
as $m\to\infty$. Thus there are unit vectors $z_m^{ji}\in\R^n$ and constants $c_m^{ji}\in\R$ as $j=1,\ldots,k$ and $i=1,\ldots,
2^m$ satisfying \eqref{3.1a} and the relationships in \eqref{3.1b}, \eqref{3.1c} with the sets $C_m^i$ defined in \eqref{3.1C}. $\h$\vspace*{-0.2in}

\section{Discrete Approximations of Local Optimal Solutions}\label{dis-app-opt}
\setcounter{equation}{0}\vspace*{-0.1in}

As seen above, Theorem~\ref{Thm3.1} provides a constructive discrete approximation of {\em any} feasible solution to problem $(P)$ by
feasible solutions to discrete-time problems, with no connections to optimization. The main goal here is to study a {\em given} local optimal
solution to $(P)$ by using discrete approximations as a {\em vehicle} to derive further necessary optimality conditions for it. To proceed in this
direction,  we construct a sequence of discrete-time optimization problems such that their optimal solutions always exist and strongly converge
in the sense below to the given local minimizer of the original sweeping control problem.

Our main attention in this section is paid to {\em relaxed $W^{1,2}\times L^2$-local minimizers} $(\ox(\cdot),\ou(\cdot))$ for $(P)$ introduced in
Definition~\ref{relaxed} while recalling that the relaxation is not needed if either the set $g(x,U)$ is convex, or $(\ox(\cdot),\ou(\cdot))$ is a
{\em strong} local minimizer for $(P)$; see the discussions in Section~\ref{stand}.

Given a relaxed $W^{1,2}\times L^2$-local minimizer $(\ox(\cdot),\ou(\cdot))$, we construct the following family of discrete-time problems $
(P_m)$, $m\in\N$, where $F$ is defined in \eqref{F0}, and where $z_m^{ji},c_m^{ji}$ are taken from Theorem~{\rm\ref{Thm3.1}}:
\begin{equation}\label{e:3.15}
\textrm{minimize }\;J_m[x_m,u_m]:=\vph\big(x_m(T)\big)+
\frac{1}{2}\sum_{i=0}^{2^m-1}\int_{t^i_m}^{t^{i+1}_m}\(\n\frac{x^{i+1}_m-x^i_m}{h_m}-
\dot{\ox}(t)\en^2+\n u^{i}_m-\ou^{i}\en^2\)dt
\end{equation}
over discrete trajectories $(x_m,u_m)=(x^0_m,x^1_m,\ldots,x^{2^m}_m,u^0_m,u^1_m,\ldots,u^{2^m-1}_m)$ subject to the constraints
\begin{equation}\label{e:3.16}
x^{i+1}_m\in x^i_m-h_m F(x^i_m,u^i_m)\textrm{ for }\;i=0,\ldots,2^m-1,
\end{equation}
\begin{equation*}
\langle z_m^{ji},x^{i}_m\rangle\le c_m^{ji}\;\textrm{ for  all }\;j=1,\ldots,s\;\textrm{ and }\;i=1,\ldots,2^m\;\text{ with }\;x^0_m:=x_0\in C,
\;u^0_m:=\ou (0),
\end{equation*}
\begin{equation}\label{e:3.18}
\sum_{i=0}^{2^m-1}\int_{t^i_m}^{t^{i+1}_m}\(\n\frac{x^{i+1}_m-x^i_m}{h_m}-\dot{\ox}(t)\en^2+\n u_m^i-\ou^i\en^2\)dt\le\frac{\e}{2},
\end{equation}
\begin{equation}\label{e:3.18+}
u^i_m\in U\;\textrm{ for }\;i=0,\ldots,2^m-1.
\end{equation}

To implement the method of discrete approximation, we have to make sure that each problem $(P_m)$ admits an optimal solution. By taking
into account Theorem~\ref{Thm3.1}, we deduce it from the classical Weierstrass existence theorem in finite dimensions due to the construction
of $(P_m)$ and the assumptions made.\vspace*{-0.07in}

\begin{proposition}\label{ThmExis2} In addition to the assumptions of Theorem~{\rm\ref{Thm3.1}}, suppose that the cost function $\ph$ is lower
semicontinuous $($l.s.c.$)$ on bounded sets. Then each problem $(P_m)$ admits an optimal solution provided that $m\in\N$ is sufficiently
large.
\end{proposition}\vspace*{-0.07in}
{\bf Proof.} It follows from Theorem~\ref{Thm3.1} that the set of feasible solutions $(x_m,u_m)$ to $(P_m)$ is nonempty for any large $m$. It
follows from the constraint structures in $(P_m)$ and the assumptions imposed on $U$ and $g$ that the feasible sets are closed. Furthermore, it
easy to deduce from the localization in \eqref{e:3.18} that the feasible sets are bounded as well. Thus the lower semicontinuity assumption on
the cost function $\ph$ ensures the existence of optimal solutions to $(P_m)$ by the Weierstrass theorem. $\h$\vspace*{0.02in}

Now we are ready to derive the main result of this section that establishes the strong $W^{1,2}$ convergence of any sequence $(\ox_m(\cdot),
\ou_m(\cdot))$ of optimal solutions to $(P_m)$, which are extended to the entire interval $[0,T]$, to the given local minimizer $(\ox(\cdot),\ou
(\cdot))$ for the original problem $(P)$.\vspace*{-0.07in}

\begin{theorem}\label{ThmStrong} Let $(\ox(\cdot),\ou(\cdot))$ be a relaxed $W^{1,2}\times W^{1,2}$-local minimizer for the sweeping control
problem $(P)$, and let $\ph$ be continuous around $\ox(T)$ in addition to the assumptions of Theorem~{\rm\ref{Thm3.1}}. Consider any
sequence of optimal solutions $(\ox_m(\cdot),\ou_m(\cdot))$ to problems $(P_m)$ and extend them to $[0,T]$ piecewise linearly for
$\ox_m (\cdot)$ and piecewise constantly for $\ou_m(\cdot)$ without relabeling. Then we have the convergence
\begin{equation*}
\big(\ox_m(\cdot),\ou_m(\cdot)\big)\to\big(\ox(\cdot),\ou(\cdot)\big)\;\textrm{ as }\;m\to\infty
\end{equation*}
in the strong topology of $W^{1,2}([0,T];\R^n)\times L^2([0,T];\R^d)$.
\end{theorem}\vspace*{-0.07in}
{\bf Proof}.
It is sufficient to show that
\begin{equation}\label{e:6.10}
\underset{m\to\infty}{\mathrm{lim}}\int_0^T\(\n\dot{\ox}_m(t)-\dot{\ox}(t)\en^2+\n\ou_m(t)-\ou(t)\en^2\)dt=0.
\end{equation}
Arguing by contradiction, suppose that there exists a subsequence of the integral values $\gg_m$ in \eqref{e:6.10} that converges, without
relabeling, to some number $\gg>0$. Due to \eqref{e:3.18}, the sequence of extended optimal solutions $\{(\dot\ox_m(\cdot),\ou_m(\cdot))\}$ to
$(P_m)$ is bounded in the reflexive space $L^2([0,T];\R^n)\times L^2([0,T];\R^d)$, and thus it contains a weakly convergence subsequence in
this product space, again without relabeling. Denote by $(\Tilde v(\cdot),\Tilde u(\cdot))$ the limit of the latter subsequence and then let
\begin{equation*}
\tilde{x}(t):=x_0+\int_0^T\tilde v(\tau)d\tau\;\textrm{ for all }\;t\in[0,T].
\end{equation*}
Since $\dot{\tilde{x}}(t)=\tilde v(t)$ for a.e.\ $t\in[0,T]$, we have that
\begin{equation*}
\big(\ox_m(\cdot),\ou_m(\cdot)\big)\to\big(\tilde x(\cdot),\tilde u(\cdot)\big)\;\textrm{ as }\;m\to\infty
\end{equation*}
in the topology of $W^{1,2}([0,T];\R^n)\times L^2([0,T];\R^d)$. Invoking the Mazur weak closure theorem tells us that there is a sequence
of convex combinations of $(\ox_m(\cdot),\ou_m(\cdot))$, which converges to $(\tilde x(\cdot),\tilde u(\cdot))$ strongly in $W^{1,2}([0,T];\R^n)
\times L^2([0,T];\R^d)$, and thus $(\dot{\ox}_m(t),\ou_m(t))\to(\dot{\tilde x}(t),\tilde u(t))$ for a.e.\ $t\in[0,T]$ along a subsequence. Furthermore,
we can clearly replace above the piecewise linear extensions of the discrete trajectories $\ox_m(\cdot)$ to the interval $[0,T]$ by the trajectories
of \eqref{Problem} generated by the controls $\ou_m(\cdot)$ piecewise constantly extended to $[0,T]$. The obtained pointwise convergence of
convex combinations allows us to conclude that $\tilde u(t)\in\co U$ for a.e.\ $t\in[0,T]$ and that $\tilde x(\cdot)$ satisfies the convexified
differential inclusion \eqref{conv}. Passing now to the limit as $m\to\infty$ in the cost functional and constraints \eqref{e:3.15}--\eqref{e:3.18+} of problem $(P_m)$ with taking into account the assumed local continuity of $\ph$ and the constructions above, we conclude that the pair
$(\tilde x (\cdot),\tilde u(\cdot))$ belongs to the prescribed $W^{1,2}\times L^2$-neighborhood of the given local minimizer
$(\ox(\cdot),\ou(\cdot))$ and satisfies the inequality
\begin{equation}\label{contr}
J[\tilde{x},\tilde{u}]+\gg/2\le J[\ox,\ou]\Longrightarrow J[\tilde{x},\tilde{u}]<J[\ox,\ou]
\end{equation}
due the aforementioned strong convergence of $(\ox_m(\cdot),\ou_m(\cdot))$ to $(\tilde x(\cdot),\tilde u(\cdot))$ and the structure of \eqref{e:3.15}. Appealing to Definition~\ref{relaxed} tells us that \eqref{contr} contradicts the very fact that $(\ox(\cdot),\ou(\cdot))$ is a relaxed $W^{1,2}\times L^2$-local minimizer of $(P)$. Thus we get \eqref{e:6.10} and complete the proof of the theorem. $\h$\vspace*{0.01in}

Recalling the discussion after Definition~\ref{relaxed} leads us to the following consequence of Theorem~\ref{ThmStrong}, which provides the
strong approximation of local minimizers for $(P)$ without an explicit relaxation.\vspace*{-0.07in}

\begin{corollary}\label{ThmStrong1} In addition to the assumptions of Theorem~{\rm\ref{ThmStrong}}, suppose that the sets $g(x,U)$ and $U$
are convex. Then the convergence result of Theorem~{\rm\ref{ThmStrong}} holds true.
\end{corollary}\vspace*{-0.3in}

\section{Tools of Variational Analysis}\label{gen-diff}
\setcounter{equation}{0}\vspace*{-0.1in}

The results of Section~\ref{dis-app-opt} make a bridge between the given local minimizer $(\ox(\cdot),\ou(\cdot))$ of the original problem $(P)$
and (global) optimal solutions for the sequence of discrete approximations $(P_m)$ that exist by Proposition~\ref{ThmExis2} and strongly
converge to $(\ox(\cdot),\ou(\cdot))$ by Theorem~\ref{ThmStrong}. This supports our approach to derive necessary optimality conditions for
$(\ox(\cdot),\ou(\cdot))$ by establishing firstly necessary conditions for optimal solutions to the discrete-time problems $(P_m)$
and then passing to the limit in them as $m\to\infty$.

Looking at the structures of each problem $(P_m)$ and the equivalent problem of finite-dimensional mathematical programming defined in
Section~\ref{nc-disc}, we observe that they are always {\em nonsmooth} and {\em nonconvex}, even when the initial data of $(P)$ possess
these properties. This is due to the {\em graphical} set constraints associated with the discrete-time inclusions \eqref{e:3.16} that are generated
by the normal cone mapping in \eqref{F0}.

To proceed with deriving necessary optimality conditions for $(P_m)$ and then for $(P)$ by passing to the limit, we have to employ appropriate
generalized differential constructions of variational analysis. These constructions should be {\em robust}, enjoy comprehensive {\em calculus
rules}, and such that the corresponding normal cone is {\em not too large} while being applied to--specifically--graphical sets. It does hold, in
particular, for the Clarke normal cone $\Bar N$, which is always a linear subspace of a maximum dimension for sets that are graphically
homeomorphic to graphs of Lipschitzian functions; see \cite{m-book1,rw} for more details and references. For example, we have $\Bar N((0,0);
\gph|x|)=\R^2$ for the graph of the simplest convex function on $\R$.

All the required properties are satisfied for the generalized differential constructions initiated by the second author. Elements of the first-order
theory and various applications can be found by now in many books; see, e.g., \cite{m-book1}--\cite{m18}, \cite{rw}, \cite{v}.  We refer the reader
to \cite{m-book2,m18} and the bibliographies therein for second-order constructions used in what follows.

To briefly overview the needed notions, recall first the (Painlev\'e-Kuratowski) {\em outer limit} of a set-valued mapping/multifunction $F\colon
\R^n \tto\R^m$ at $\ox$ with $F(\ox)\ne\emp$ given by
\begin{equation}\label{c53}
\underset{x\to\ox}{\textrm{Lim sup }}F(x):=\big\{y\in\R^m\;\big|\;\exists\textrm{ sequences }\;x_k\to\ox,\;y_k\to y\;\textrm{ such that }\;y_k\in
F(x_k),\;k \in\N\big\}
\end{equation}
Given now a set $\O\subset\R^n$ locally closed around $\ox\in\O$, we define by using \eqref{c53} the (basic, limiting, Mordukhovich) {\em
normal cone} to $\O$ at $\ox$ by
\begin{equation}\label{c54}
N(\ox;\O)=N_\O(\ox):=\underset{x\to\ox}{\textrm{Lim sup}}\big\{\textrm{cone}[x-\Pi(x;\O)]\big\}.
\end{equation}
where $\Pi(x;\O):=\big\{u\in\O\;\big|\;\|x-u\|=\dist(x;\O)\big\}$ is the Euclidean projection of $x$ onto $\O$, and where `cone' stands for the
(nonconvex) conic hull of the set. When $\O$ is convex, $(\ref{c54})$ reduces to the normal cone of convex analysis, but it is often nonconvex
otherwise.

Given further a set-valued mapping $F\colon\R^n\tto\R^m$ with its domain and graph
\begin{equation*}
\dom F:=\big\{x\in\R^n\;\big|\;F(x)\ne\emp\}\;\textrm{ and }\;\gph F:=\big\{(x,y)\in\R^n\times\R^m\;\big|\;y\in F(x)\big\}
\end{equation*}
locally closed around $(\ox,\oy)\in\gph F$, the {\em coderivative} of $F$ at $(\ox,\oy)$ is generated by \eqref{c54} as
\begin{equation}\label{c55}
D^*F(\ox,\oy)(u):=\big\{v\in\R^n\;\big|\;(v,-u)\in N\big((\ox,\oy);\gph F\big)\big\},\quad u\in\R^m.
\end{equation}
When $F\colon\R^n\to\R^m$ is single-valued and continuously differentiable $({\cal C}^1$-smooth) around $\ox$, we have
\begin{equation*}
D^*F(\ox)(u)=\big\{\nabla F(\ox)^*u \big\}\;\textrm{ for all }\;u\in\R^m
\end{equation*}
via the adjoint/transposed Jacobian matrix $\nabla F(\ox)^*$, where $\oy=F(\ox)$ is omitted.

Let $\phi\colon\R^n\to\oR:=(-\infty,\infty]$ be an extended-real-valued l.s.c.\ function $\phi\colon\R^n\to\oR:=(-\infty,\infty]$ with
\begin{equation*}
\dom\phi:=\big\{x\in\R^n\;\big|\;\vph(x)<\infty\big\}\;\textrm{ and }\;\epi\phi:=\big\{(x,\alpha)\in\R^{n+1}\;\big|\;\alpha\ge\phi(x)\big\}
\end{equation*}
standing for its domain and epigraph. The (first-order) {\em subdifferential} of $\phi$ at $\ox\in\dom\phi$ is defined geometrically
via the normal cone \eqref{c54} by
\begin{equation}\label{sub1}
\partial\phi(\ox):=\big\{v\in\R^m\;\big|\;(v,-1)\in N\big((\ox,\phi(\ox));\epi\phi\big)\big\}
\end{equation}
while admitting equivalent analytic representations; see, e.g., \cite{m-book1,rw}. Note that $N(\ox;\O)=\partial\dd(\ox;\O)$ for any $\ox\in\O$,
where $\dd(x;\O)$ denotes the indicator function of $\O$ equal to 0 for $x\in\O$ and $\infty$ otherwise. Then given a
subgradient $\ov\in\partial \phi(\ox)$ and following \cite{m-book1,m18}, we define the {\em second-order subdifferential}
(or {\em generalized Hessian}) of $\phi$ at $\ox$ relative to $\ov$ by
\begin{equation*}
\partial^2\phi(\ox,\ov)(u):=(D^*\partial\phi)(\ox,\ov)(u),\quad u\in\R^n,
\end{equation*}
via the coderivative \eqref{c55} of the first-order subdifferential mapping $x\mapsto\partial\phi(x)$ from \eqref{sub1}. If the function $\phi$ is
${\cal C}^2$-smooth around $\ox$, then we have the representation
\begin{equation*}
\partial^2\phi(\ox,\ov)(u)=\big\{\nabla^2\phi(\ox)u\big\}\;\textrm{ for all }\;u\in\R^n,
\end{equation*}
where $\nabla^2\phi(\ox)$ stands for the classical (symmetric) Hessian of $\phi$ at $\ox$ with $\ov=\nabla\phi(\ox)$. If $\phi(x):=\dd(x;\O)$, then
$\partial^2\ph(\ox,\ov)(u)=(D^*N_\O)(\ox,\ov)(u)$ for any $\ov\in N(\ox;\O)$ and $u\in\R^n$. The latter second-order construction is evaluated
below in the case of the polyhedral set $\O=C$ from \eqref{C}. To proceed, define the index sets corresponding the generating vectors $x^j_*$
in \eqref{C} by
\begin{equation}\label{c56}
I_0(w):=\big\{j\in I(x)\;\big|\;\la x^j_*,w\ra=c_j\big\}\;\textrm{ and }\;I_>(w):=\big\{j\in I(x)\;\big|\;\la x^j_*,w\ra>c_j\big\},\;w\in\R^n.
\end{equation}
where $I(x)$ is taken from \eqref{aci} with $\ox:=x\in C$. The next theorem provides an effective upper estimate of the coderivative of $F$ from
\eqref{F0} with ensuring the equality under an additional assumption on $x^j_*$.\vspace*{-0.07in}

\begin{theorem}\label{Thm6.1}
Given $F$ in $(\ref{F0})$ with $C$ from $(\ref{C})$, denote $G(x):= N(x;C)$ and suppose in addition to standing assumptions that $g$ is
${\cal C}^1$-smooth around the reference points. Then for any $(x,u)\in C\times U$ and $\o+g(x,u)\in G(x)$ we have the
coderivative upper estimate
\begin{equation}\label{c57}
D^*F(x,u,\o)(w)\subset\Big\{z=\Big(-\nabla_x g(x,u)^*w+\underset{j\in I_0(w)\cup I_>(w)}\sum\gg^j x^j_*,-\nabla_u g(x,u)^*w\Big)\Big\},
\end{equation}
where $w\in\dom D^*G(x,\o+g(x,u))$, where $I_0(w)$ and $I_>(w)$ are taken from \eqref{c56}, and where $\gg^j\in\R$ for $j\in I_0(w)$,
while $\gg^j\ge 0$ for $j\in I_>(w)$. Furthermore, \eqref{c57} holds as an equality and the domain $\dom D^*G(x,\o+g(x,u))$
can be computed by
\begin{equation}\label{c58}
\dom D^*G\big(x,\o+g(x,u)\big)=\Big\{w\Big|\;\exists\lm^j\ge 0\textrm{ with }\o+g(x,u)=\underset{j\in I(x)}\sum\lm^j x^j_*,\;\lm^j>0
\Longrightarrow\la x^j_*,w\ra=c_j\Big\}
\end{equation}
provided that the generating vectors $\{x^j_*\;|\;j\in I(x)\}$ of the polyhedron $C$ are linearly independent.
\end{theorem}\vspace*{-0.07in}
{\bf Proof}. Picking any $w\in\dom D^*G(x,\o+g(x,u))$ and $z\in D^*F(x,u,y)(w)$ and then denoting $\tilde{G}(x,u):=G(x)$ and
$\tilde{f}(x,u):=-g(x,u)$, we deduce from \cite[Theorem~1.62]{m-book2} that
\begin{equation*}
z\in\nabla\tilde{f}(x,u)^*w+D^*\tilde{G}\big(x,u,\o+g(x,u)\big)(w).
\end{equation*}
Observe then the obvious composition representation
\begin{equation*}
\tilde{G}(x,u)=G\circ\tilde{g}(x,u)\;\textrm{ with }\;\tilde{g}(x,u):=x,
\end{equation*}
where the latter mapping has the surjective derivative. It follows from \cite[Theorem~1.66]{m-book2} that
\begin{equation}\label{59}
z\in\nabla\tilde{f}(x,u)^*w+\nabla\tilde{g}(x,u)^*D^*G(\big(x,\o+g(x,u)\big)(w).
\end{equation}
Employing now in \eqref{59} the coderivative estimate for the normal cone mapping $G$ obtained in \cite[Theorem~4.5]{hmn} with the exact
coderivative calculation given in \cite[Theorem~4.6]{hmn} under the linear independence of the generating vectors $x^j_*$ and also taking into
account the structure of the mapping $\tilde{f}$ in \eqref{59}, we arrive at \eqref{c57} and the equality therein under the aforementioned
assumption. $\h$\vspace*{-0.2in}

\section{Necessary Optimality Conditions for Discrete-Time Problems}\label{nc-disc}
\setcounter{equation}{0}\vspace*{-0.1in}

Here we derive necessary optimality conditions for solutions to each problem $(P_m)$, $m\in\N$, formulated in \eqref{e:3.15}--\eqref{e:3.18+}. It
will be done by reducing each $(P_m)$ to a nondynamic problem of nondifferentiable programming with functional and many geometric
constraints, then employing necessary optimality conditions for the latter problem obtained in terms of generalized differential constructions of
Section~\ref{gen-diff}, and finally expressing the obtained conditions in terms of the given data of $(P_m)$ by using calculus rules of
generalized differentiation. In this way we arrive at the following necessary conditions, which will be further specified below by applying the
second-order calculations presented in Section~\ref{gen-diff}.\vspace*{-0.07in}

\begin{theorem}\label{Thm5.2*}
Let $(\ox_m,\ou_m)=(\ox^0_m,\ldots,\ox^{2^m}_m,\ou^0_m,\ldots,\ou^{2^m-1}_m)$ be an optimal solution to problem $(P_m)$. Assume that
$\gph F$ is closed and the function $\vph$ is Lipschitz continuous around the point $\ox_m(T)$. Then there are elements $\lm_m\ge 0$,
$\psi_m=(\psi^0_m,\ldots,\psi^{2^m-1}_m)$ with $\psi^i_m\in N(\ou^i_m;U)$ as $i=0,\ldots,2^m-1,$ $\xi_m=(\xi^{1}_m,\ldots,
\xi^{s}_m)\in\R^s_+$,  and $p^i_m\in\R^n $ as $i=0,\ldots,2^m$ satisfying the conditions
\begin{equation}\label{e:5.8*}
\lm_m+\n\xi_m\en+\sum_{i=0}^{2^m-1}\n p^{i}_m\en+\n\psi_m\en\ne 0,
\end{equation}
\begin{equation}\label{xi}
\xi^{j}_m\big(\la z^{j2^m}_m,x^{2^m}_m\ra-c^{j2^m}_m\big)=0,\quad j=1,\ldots,s,
\end{equation}
\begin{equation}\label{mutx}
-p^{2^m}_m=\lm_m\vth^{2^m}_m+\sum_{j=1}^s\xi^{j}_m z^{j2^m}_m\in\lm_m\partial\vph(\ox^{2^m}_m)+\sum_{j=1}^s\xi^{j}_m z^{j2^m}_m,
\end{equation}
\begin{equation}\label{e:5.10*}
\begin{array}{ll}
&\disp\Big(\frac{p^{i+1}_m-p^{i}_m}{h_m},-\frac{1}{h_m}\lm_m\th^{iu}_m,\frac{1}{h_m}\lm_m\th^{iy}_m-p^{i+1}_m\Big)\\
&\in\disp\Big(0,\frac{1}{h_m}\psi^i_m,0\Big)+N\Big(\Big(\ox^i_m,\ou^i_m,-\frac{\ox^{i+1}_m-\ox^i_m}{h_m}\Big);\gph F\Big)
\end{array}
\end{equation}
for $i=0,\ldots,2^m-1$, where we use the notation
\begin{equation}\label{theta}
\theta^i_m=\(\th^{iy}_m,\th^{iu}_m\):=\Big(\int_{t^i_{m}}^{t^{i+1}_{m}}\Big\|\frac{\ox^{i+1}_m-\ox^i_m}{h_m}-\dot{\ox}(t)\Big\| dt,\int_{t^i_{m}}^{t^{i+1}_{m}}
\|\ou^i_m-\ou(t)\|dt\Big).
\end{equation}
\end{theorem}\vspace*{-0.07in}
{\bf Proof.} Denote $z:=(x^0_m,\ldots,x^{2^m}_m,u^0_m,\ldots,u^{2^m-1}_m,y^0_m,\ldots,y^{2^m-1}_m)\in\R^{(2\cdot 2^m+1)n+2^m\cdot d}$,
where the starting point $x^0_m$ is fixed. Taking $\e>0$ from $(P_m)$, consider the following problem of mathematical programming $(MP)$
with respect to the variable $z$:
\begin{equation*}
\textrm{minimize }\;\phi_0(z):=\vph\big(x(T)\big)+\frac{1}{2}\sum_{i=0}^{2^m-1}\int_{t^i_m}^{t^{i+1}_m}\n\(y^i_m-\dot{\ox}(t),u_m^i-\ou(t)\)\en^2dt
\end{equation*}
subject to finitely many equality, inequality, and geometric constraints given by
\begin{equation*}
\phi(z):=\sum_{i=0}^{2^m-1}\int_{t^i_m}^{t^{i+1}_m}\n\(y^i_m,u^i_m\)-\big(\dot{\ox}(t),\ou(t)\big)\en^2dt-\frac{\e}{2}\le 0,
\end{equation*}
\begin{equation*}
g_i(z):=x^{i+1}_m-x^i_m-h_m y^i_m=0,\quad i=0,\ldots,2^m-1,
\end{equation*}
\begin{equation*}
h_{j}(z):=\la z^{j2^m}_m,x^{2^m}_m\ra-c^{j2^m}_m \le 0,\quad j=1,\ldots,s,
\end{equation*}
\begin{equation*}
z\in\Xi_i:=\left\{(x^0_m,\ldots,y^{2^m-1}_m)\in\R^{(2\cdot 2^m+1)n+2^m\cdot d}\;\Big|\;-y^i_m
\in F(x^i_m,u^i_m)\right\},\quad i=0,\ldots,2^m-1,
\end{equation*}
\begin{equation*}
z\in\Xi_{2^m}:=\big\{(x^0_m,\ldots,y^{2^m-1}_m)\in\R^{(2\cdot 2^m+1)n+2^m\cdot d}\;\big|\;x^0_m\;\textrm{ is fixed}\big\},
\end{equation*}
\begin{equation*}
z\in\O_i:=\big\{(x^0_m,\ldots,y^{2^m-1}_m)\in\R^{(2\cdot 2^m+1)n+2^m\cdot d}\;\big|\;u^i_m\in U\big\},\quad i=0,\ldots,2^m-1,
\end{equation*}
Necessary optimality conditions for problem $(MP)$ in terms of the generalized differential tools reviewed above can deduced from
\cite[Proposition~6.4 and Theorem~6.5]{m18}. We specify them for the optimal solution
$$
\oz:=\(\ox^0_m,\ldots,\ox^{2^m}_m,\ou^0_m,\ldots,\ou^{2^m-1}_m,\oy^0_m,\ldots,\oy^{2^m-1}_m\)
$$
to $(MP)$. It follows from Theorem~\ref{ThmStrong} that the inequality constraint in $(MP)$ defined by $\phi$ is inactive for large $m$, and so
the corresponding multiplier does not appear in the optimality conditions. Thus we can find $\lm_m\ge 0$, $\xi_m=(\xi^{1}_m,\ldots,\xi^{s}_m)\in
\R^s_+$, $p^i_{m}\in\R^{n}$ as $i=1,\ldots,2^m$, and
$$
z^*_i=\big(x^*_{0i},\ldots,x^*_{2^mi},u^*_{0i},\ldots,u^*_{(2^m-1)i},y^*_{0i},y^*_{1i},\ldots,y^*_{(2^m-1)i}\big),\quad i=0,\ldots,2^m,
$$
which are not zero simultaneously while satisfying the conditions
\begin{equation}\label{69}
z^*_i\in\left\{\begin{matrix}
N(\oz;\Xi_i)+N(\oz;\O_i)\;\textrm{ if }\;i\in\big\{0,\ldots,2^m-1\big\},\\
N(\oz;\Xi_i);\textrm{ if }\;i=2^m,
\end{matrix}\right.
\end{equation}
\begin{equation*}
-z^*_0-\ldots-z^*_{2^m}\in\lm_m\partial\phi_0(\oz)+\sum_{j=1}^{s} \xi^{j}_m\nabla h_{j}(\oz)+\sum_{i=0}^{2^m-1}\nabla g_i(\oz)^*p^{i+1}_m,
\end{equation*}
\begin{equation}\label{71+}
\xi^{j}_m h_{j}(\oz)=0,\quad j=1,\ldots,s.
\end{equation}
Note that the first line in \eqref{69} comes by applying the normal cone intersection formula from \cite[Corollary~3.5]{m-book1} to
$\oz\in\O_i\cap \Xi_i$ for $i=0,\ldots,2^m-1$. It follows from the structure of the sets $\O_i$ and $\Xi_i$ that the inclusions in \eqref{69}
can be equivalently written as
\begin{equation}\label{e:5.18*}
\big(x^*_{ii},u^*_{ii}-\psi^{i}_m,-y^*_{ii}\big)\in N\Big(\Big(\ox^i_m,\ou^i_m,-\frac{\ox^{i+1}_m-\ox^i_m}{h_m}\Big);\gph F\Big)\;\textrm{ for }
\;i=0,\ldots,2^m-1
\end{equation}
with every other components of $z^*_i$ equal to zero, where $\psi^{i}_m\in N(\ou^i_m;U)$ for all $i=0,\ldots,2^m-1$. Observe furthermore that
$x^*_{0m}$ and $u^*_{0m}$ determined by the normal cone to $\Xi_m$ are the only nonzero components of $z^*_m$. This implies by using
\eqref{69}--\eqref{71+} that
$$
-z^*_0-\ldots-z^*_{2^m}\in\lm_m\partial\phi_0(\oz)+\sum_{j=1}^{s}\xi^{j}_m\nabla h_{j}(\oz)+\sum_{i=0}^{2^m-1}\nabla g_i(\oz)^*p^{i+1}_m,
$$
with $\xi^{j}_m\(\la z^{j2^m}_m, x^{2^m}_m\ra-c^{j2^m}_m\)=0,\;j=1,\ldots,s$. Using the expressions for $\phi_0$, $g_i$, and $h_j$ above
together with the elementary subdifferential sum rule from \cite[Proposition~1.107]{m-book1} gives the calculations
$$
\Big(\sum_{j=1}^{s}\xi^{j}_m\nabla h_{j}(\oz)\Big)_{x^{2^m}_m}=\Big(\sum_{j=1}^{s}\xi^{j}_m z^{j2^m}_m\Big),
$$
$$
\Big(\sum_{i=0}^{2^m-1}\nabla g_i(\oz)^*p^{i+1}_m\Big)_{x^i_m}=\left\{\begin{matrix}
-p^{1}_m;\textrm{ if }\;i=0,\\
p^{i}_m-p^{i+1}_m\;\textrm{ if }\;i=1,\ldots,2^m-1,\\
p^{2^m}_m\;\textrm{ if }\;i=2^m,
\end{matrix}\right.
$$
$$
\Big(\sum_{i=0}^{2^m-1}\nabla g_i(\oz)^*p^{i+1}_m\Big)_{y^i_m}=\(-h_m p^{1}_m,-h_m p^{2}_m,\ldots,-h_m p^{2^m}_m\),
$$
$$
\partial(\phi_0)(\oz)=\partial\vph(\ox^m_m)+\frac{1}{2}\sum_{i=0}^{2^m-1}\nabla\rho_i(\oz)\;\textrm{ with }\;\rho_i(\oz):=\int_{t^i_m}^{t^{i+1}_m}
\Big\| \Big(\frac{\ox^{i+1}_m-\ox^i_m}{h_m}-\dot{\ox}(t),\ou^i_m-\ou(t)\Big)\Big\|^2dt.
$$
The set $\lm_m\partial\phi_0(\oz)$ is represented as the collection of
$$
\lm_m\big(0,\ldots,0,\vth^{2^m}_m,\th^{0u}_m,\ldots,\th^{(2^m-1)u}_m,\th^{0y}_m,\ldots,\th^{(2^m-1)y}_{m}\big)\;\textrm{ with }\;\vth^{2^m}_m\in
\partial\vph(\ox^{2^m}_m),
$$
$$
(\th^{iu}_m,\th^{iy}_m)=\Big(\int_{t^i_m}^{t^{i+1}_m}\|\ou^i_m -\ou(t)\| dt,\;\int_{t^i_m}^{t^{i+1}_m}\Big\|\frac{\ox^{i+1}_m-\ox^i_m}{h_m}-\dot{\ox}(t)
\Big\|dt\Big),\quad i=0,\ldots,2^m-1.
$$
Thus we obtain the following relationships
\begin{equation}\label{daux}
-x^*_{00}-x^*_{02^m}=-p^{1x}_m,
\end{equation}
\begin{equation}\label{e:5.21*}
-x^*_{ii}=p^{i}_m-p^{i+1}_m,\quad i=1,\ldots,2^m-1,
\end{equation}
\begin{equation}\label{e:5.24*}
0=\lm_m\vth^{2^m}_m+p^{2^m}_m+\sum_{j=1}^{s}\xi^{j}_m z^{j2^m}_m\;\textrm{ with }\;\vth^{2^m}_m\in\partial\vph(\ox_m^{2^m}),
\end{equation}
\begin{equation}\label{e:5.22*}
-u^*_{00}=\lm_m\th^{0u}_m\;\textrm{ and }\;-u^*_{ii}=\lm_m\th^{iu}_m,\quad i=1,\ldots,2^m-1,
\end{equation}
\begin{equation}\label{e:5.23*}
-y^*_{ii}=\lm_m\th^{iy}_m-h_m p^{i+1}_m,\quad i=0,\ldots,2^m-1,
\end{equation}
which allow us to arrive at all the necessary optimality conditions claimed in the theorem. Indeed, observe first that \eqref{71+} yields \eqref{xi}.
Extending $p_m$ by $p^0_m:=x^*_{02^m}$ ensures that \eqref{mutx} follows from $(\ref{e:5.24*})$.
Then we deduce from \eqref{e:5.21*}, \eqref{e:5.22*}, and \eqref{e:5.23*} that
$$
\frac{x^*_{ii}}{h_m}=\frac{p^{i+1}_m-p^{i}_m}{h_m},\;\;\frac{u^*_{ii}}{h_m}=-\frac{1}{h_m}\lm_m\th^{iu}_m,\;\textrm{ and }\;\frac{y^*_{ii}}{h_m}
=-\frac{1}{h_m}\lm_m\th^{iy}_m+p^{i+1}_m.
$$
Substituting this into the left-hand side of \eqref{e:5.18*} justifies the discrete-time adjoint inclusion \eqref{e:5.10*}.

Finally, to verify \eqref{e:5.8*} we argue by contradiction and suppose that $\lm_m=0,\xi_m=0,\psi_m=0$, and $p^{i}_m=0$ as
$i=0,\ldots, 2^m-1$, which yield $x^*_{02^m}=p^{0}_m=0$. Then it follows from \eqref{e:5.24*} that $p^{2^m}_m=0$, and so
$p^{i}_m=0$ whenever
$i=0,\ldots,2^m$. By \eqref{daux} and \eqref{e:5.21*} we get $x^*_{ii}=0$ for all $i=0,\ldots,2^m-1$.
Using \eqref{e:5.22*} tells us that $u^*_{ii} =0$ as $i=1,\ldots,2^m-1$. Since the first condition in \eqref{e:5.22*} yields also $u^*_{00}=0$,
it follows that $u^*_{ii}=0$ for $i=0,\ldots,2^m-1$.
In addition we have by \eqref{e:5.23*} that $y^*_{ii}=0$ for all $i=0,\ldots,2^m-1$. Remembering that the components of $z^*_i$
different from $(x^*_{ii},u^*_{ii},y^*_{ii})$ are zero for $i=0,\ldots,2^m-1$ ensures that $z^*_{i}=0$ for $i=0,\ldots,2^m-1$
and similarly $z^*_{2^m}=0$.
Therefore $z^*_i=0$ for all $i=0,\ldots,2^m$, which violates the nontriviality condition for $(MP)$ and thus completes the
proof. $\h$\vspace*{0.01in}

The next theorem applies to \eqref{e:5.10*} the calculation result of Theorem~\ref{ThmStrong} and provides in this way necessary optimality
conditions for problem $(P_m)$ expressed entirely via its initial data.\vspace*{-0.1in}

\begin{theorem}\label{Thm6.2} Let $(\ox_m,\ou_m)$ be an optimal solution to problem $(P_m)$ formulated in \eqref{e:3.15}--\eqref{e:3.18+},
where the cost function $\ph$ is locally Lipschitzian around $\ox_m(T)$, and where the sweeping mapping $F$ is defined in \eqref{F0}. Using
the notation and assumptions of Theorem~{\rm\ref{Thm6.1}}, take $(\th^{iu}_m,\th^{iy}_m)$ from \eqref{theta}. Then there exists dual elements
$(\lm_m,\psi_m,p_m)$ as in Theorem~{\rm\ref{Thm5.2*}} together with vectors $\eta^{i}_m\in\R^s_+$ for $i=0,\ldots,2^m$ and $\gg^i_m\in\R^s$
for $i=0,\ldots,2^m-1$ satisfying the nontriviality condition
\begin{equation}\label{ntc}
\lm_m+\n\eta^{2^m}_m\en+\sum_{i=0}^{2^m-1}\n p^{i}_m\en+\n\psi_m\en\ne 0,
\end{equation}
the primal-dual relationships given for all $i=0,\ldots,2^m-1$ and $j=1,\ldots,s$ by
\begin{equation}\label{87}
-\frac{\ox^{i+1}_m-\ox^i_m}{h_m}+g(\ox^i_m,\ou^i_m)=\sum_{j\in I(\ox^i_m)}\eta^{ij}_m z^{ji}_m,
\end{equation}
\begin{equation}\label{conx}
\begin{array}{ll}
\disp\frac{p^{i+1}_m-p^{i}_m}{h_m}&=-\nabla_x g(\ox^i_m,\ou^i_m)^*\Big(-\disp\frac{1}{h_m}\lm_m\th^{iy}_m+p^{i+1}_m\Big)\\
&+\disp\sum_{j\in I_0(p^{i+1}_m-\frac{1}{h_m}\lm_m\th^{iy}_m)\cup I_>(p^{i+1}_m-\frac{1}{h_m}\lm_m\th^{iy}_m)}\gg^{ij}_m z^{ji}_m,
\end{array}
\end{equation}
\begin{equation}\label{cony}
-\frac{1}{h_m}\lm_m\th^{iu}_m-\frac{1}{h_m}\psi^i_m=-\nabla_u g(\ox^i_m,\ou^i_m)^*\Big(-\frac{1}{h_m}\lm_m\th^{iy}_m+p^{i+1}_m\Big)
\end{equation}
with $\psi^{i}_m\in N(\ou^i_m;U)$ as $i=0,\ldots,2^m-1$ taken from Theorem~{\rm\ref{Thm5.2*}}, the transversality condition
\begin{equation}\label{nmutx}
-p^{2^m}_m=\lm_m\vth^{2^m}_m+\sum_{j=1}^s\eta^{2^mj}_m z^{j2^m}_m\in\lm_m\partial\vph(\ox^{2^m}_m)+\sum_{j=1}^s\eta^{2^mj}_m
z^{j2^m}_m,
\end{equation}
and such that the following implications hold for $i=0,\ldots,2^m-1$ and $j=1,\ldots,s$:
\begin{equation}\label{eta}
\Big[\la z^{ji}_m,\ox^i_m\ra<c^{ji}_m\Big]\Longrightarrow\eta^{ij}_m=0,
\end{equation}
\begin{equation}\label{93}
\left\{\begin{matrix}
\Big[j\in I_0( p^{i+1}_m-\disp\frac{1}{h_m}\lm_m\th^{iy}_m)\Big]\Longrightarrow\gg^{ij}_m\in\R,\\
\Big[j\in I_>(p^{i+1}_m-\disp\frac{1}{h_m}\lm_m \th^{iy}_m)\Big]\Longrightarrow\gg^{ij}_m\ge 0,\\
\Big[j\notin I_0(p^{i+1}_m-\disp\frac{1}{h_m}\lm_m\th^{iy}_m)\cup I_>(p^{i+1}_m-\frac{1}{h_m}\lm_m\th^{iy}_m)\Big]
\Longrightarrow\gg^{ij}_m=0.
\end{matrix}\right.
\end{equation}
We also have the complementary slackness condition $(\ref{xi})$ together with
\begin{equation}\label{94}
\[\la z^{ji}_m,\ox^i_m\ra<c^{ji}_m\]\Longrightarrow\gg^{ij}_m=0\;\textrm{ for }\;i=0,\ldots,2^m-1\;\textrm{ and }\;j=1,\ldots,s,
\end{equation}
\begin{equation}\label{eta1}
\big[\la z^{j2^m}_m,\ox^{2^m}_m\ra<c^{j2^m}_m\big]\Longrightarrow\eta^{2^mj}_m=0\;\textrm{ for }\;j=1,\ldots,s,
\end{equation}
Furthermore, the linear independence of the vectors $\{z^{ji}_m\;|\;j\in I(\ox^i_m)\}$
ensures the implication
\begin{equation}\label{96}
\eta^{ij}_m>0\Longrightarrow\Big[\Big\la z^{ji}_m,p^{i+1}_m-\disp\frac{1}{h_m}\lm_m\th^{iy}_m\Big\ra=c^{ji}_m\Big]
\end{equation}
Assuming in addition that the matrices $\nabla_u g(\ox^i_m,\ou^i_m)$ are of full rank for all $i=0,\ldots,2^m-1$ and $m\in\N$ sufficiently large,
we get the enhanced nontriviality condition
\begin{equation}\label{entc}
\lm_m+\|p^{0}_m\|+\|\psi_m\|\ne 0.
\end{equation}
\end{theorem}\vspace*{-0.07in}
{\bf Proof.} Using the necessary optimality conditions of Theorem~\ref{Thm6.1}, we can rewrite \eqref{e:5.10*} as
\begin{equation}\label{cod-disc}
\Big(\frac{p^{i+1}_m-p^{i}_m}{h_m},-\frac{1}{h_m}\lm_m\th^{iu}_m-\frac{1}{h_m}\psi^i_m\Big)\in D^*F\Big(\ox^i_m,\ou^i_m,-\frac{\ox^{i+1}_m-
\ox^i_m}{h_m}\Big) \(-\frac{1}{h_m}\lm_m\th^{iy}_m+p^{i+1}_m\)
\end{equation}
for all $i=0,\ldots,2^m-1$ by the coderivative definition \eqref{c55}. Taking into account that
\begin{equation}
-\frac{\ox^{i+1}_m-\ox^i_m}{h_m}+g(\ox^i_m,\ou^i_m)\in G(\ox^i_m)\;\textrm{ for }\;i=0,\ldots,2^m-1
\end{equation}
with $G(x)= N(x;C)$, we find vectors $\eta^{i}_m \in\R^s_+$ as $i=0,\ldots,2^m-1$ such that conditions \eqref{87} and \eqref{eta} hold.
Employing now the coderivative evaluation \eqref{c57} from Theorem~\ref{Thm6.1} with $x:=\ox^i_m$, $u:=\ou^i_m$,
$\o:=-\frac{\ox^{i+1}_m-\ox^i_m}{h_m}$, and $w:=-\frac{1}{h_m}\lm_m\th^{iy}_m+p^{i+1}_m$ for $i=0,\ldots,2^m-1$ gives us
$\gg^i_m\in\R^s$ and the relationships
$$
\Big(\frac{p^{i+1}_m-p^{i}_m}{h_m},-\frac{1}{h_m}\lm_m\th^{iu}_m-\frac{\psi^{iu}_m }{h_m}\Big)
$$
=$\(\begin{matrix}
\disp-\nabla_x g(\ox^i_m,\ou^i_m)^*\Big(-\frac{1}{h_m}\lm_m\th^{iy}_m+p^{i+1}_m\Big)+
\sum_{j\in I_0(p^{i+1}_m-\frac{1}{h_m}\lm_m\th^{iy}_m)
\cup I_>(p^{i+1}_m -\frac{1}{h_m}\lm_m\th^{iy}_m)}\gg^{ij}_m z^{ji}_m,\\
-\nabla_u g(\ox^i_m,\ou^i_m)^*\Big(-\frac{1}{h_m}\lm_m\th^{iy}_m+p^{i+1}_m\Big)
\end{matrix}
\),$\\
$$
\psi^{iu}_m\(\nu-\ou^i_m\)\le 0\;\textrm{ for all }\;\nu\in U\;\textrm{ and }\;i=0,\ldots,2^m-1.
$$
This ensures the validity of all the conditions in \eqref{conx}, \eqref{cony}, \eqref{93}, and \eqref{94}. Denoting $\eta^{2^m}_m:=\xi_m$ with
$\xi_m$ taken from Theorem~\ref{Thm5.2*}, we get $\eta^{i}_m\in\R^s_+$ for all $i=0,\ldots,2^m$ and deduce \eqref{ntc} and \eqref{nmutx} from
those in \eqref{e:5.8*} and \eqref{mutx}. Implication \eqref{eta1} follows directly from \eqref{xi} and the definition of $\eta^{2^m}_m$.

Assume finally that the generating vectors $\{z^{ji}_m\;|\;j\in I(\ox^i_m)\}$ are linear independent. In this case we deduce from \eqref{c58} and
\eqref{cod-disc} that condition \eqref{96} is satisfied. It remains to verify the enhanced nontriviality \eqref{entc} under the additional assumption
on the full rank of the matrices $\nabla_u g(\ox^i_m,\ou^i_m)$. Suppose on the contrary that $\lm_m=0$, $p^{0}_m=0$, and $\psi_m=0$. Then
$p^{i+1}_m=0$ as $i=0,\ldots,2^m-1$ by \eqref{cony}. Then it follows from \eqref{conx} the equality
\begin{equation*}
\sum_{j\in I_0(p^{i+1}_m-\frac{1}{h_m}\lm_m\th^{iy}_m)\cup I_>(p^{i+1}_m-\frac{1}{h_m}\lm_m\th^{iy}_m)}\gg^{ij}_m z^{ji}_m=0.
\end{equation*}
Invoking now \eqref{nmutx} and $p^{2^m}_m=0$ tells us that $\sum_{j=1}^s\eta^{2^mj}_m z^{j2^m}_m=0$. This implies by definition \eqref{aci}
of the active constraint indices and the imposed linear independence of $z^{ji}_m$ over this index set that $\eta^{2^m}_m=0$. Thus \eqref{ntc}
is violated, which verifies \eqref{entc} and completes the proof of the theorem. $\h$\vspace*{-0.2in}

\section{Optimality Conditions for the Controlled Sweeping Process}\label{nc-sweep}
\setcounter{equation}{0}\vspace*{-0.1in}

In this section we derive necessary optimality conditions for the local minimizer under consideration in the original problem $(P)$ by passing to
the limit as $m\to\infty$ in the necessary optimality conditions of Theorem~\ref{Thm5.2*} for the discrete-time problems $(P_m)$. Furnishing the
limiting procedure requires
the usage of Theorem~\ref{ThmStrong} and the tools of generalized differentiation reviewed in Section~\ref{gen-diff}.\vspace*{-0.1in}

\begin{theorem}\label{Thm6.1*}
Let $(\ox(\cdot),\ou(\cdot))$ be a relaxed $W^{1,2}\times L^2$-local minimizer of problem $(P)$ such that $\ou(\cdot)$ is of bounded variation
and admits a right continuous representative on $[0,T]$. In addition to $(H1)$ and $(H2)$, assume that LICQ holds along $\ox(\cdot)$ on $[0,T]
$, that $g(\cdot,\cdot)$ is ${\cal C}^1$-smooth around $(\ox(t),\ou(t))$ with the full rank of the matrices $\nabla_u g(\ox(t),\ou(t))$ on $[0,T]$, and
that $\ph$ is locally Lipschitzian around $\ox(T)$. Then there exist a multiplier $\lm\ge 0$, a signed vector measure $\gg=(\gg^1,\ldots,\gg^n)\in
C^*([0,T];\R^n)$ as well as adjoint arcs $p(\cdot)\in W^{1,2}([0,T];\R^n)$ and $q(\cdot)\in BV([0,T];\R^n)$ such that the following conditions are
fulfilled:\vspace*{-0.1in}
\begin{itemize}
\item[\bf(i)] The {\sc primal-dual dynamic relationships:}
\begin{equation}\label{37}
-\dot{\ox}(t)=\sum_{j=1}^s\eta^j(t)x^j_*-g\big(\ox(t),\ou(t)\big)\;\textrm{ for a.e. }\;t\in[0,T],
\end{equation}
where the functions $\eta^j(\cdot)\in L^2([0,T];\R_+)$ are well defined at $t=T$ while being uniquely determined by the representation in
\eqref{37};
\begin{equation}\label{c:6.6}
\dot{p}(t)=-\nabla_x g\big(\ox(t),\ou(t)\big)^*q(t)\;\textrm{ for a.e. }\;t\in[0,T],
\end{equation}
where the right continuous representative of $q(\cdot)$, with the same notation, satisfies
\begin{equation}\label{c:6.9}
q(t)=p(t)-\int_{(t,T]}d\gg(\tau)
\end{equation}
for all $t\in[0,T]$ except at most a countable subset;
\begin{equation}\label{c:6.6'}
\psi(t):=\nabla_u g\big(\ox(t),\ou(t)\big)^*q(t)\in N\big(\ou(t);U\big)\;\textrm{ for a.e. }\;t\in[0,T],
\end{equation}
which gives us the {\sc maximization condition}
\begin{equation}\label{max}
\big\la\psi(t),\ou(t)\big\ra=\max_{u\in U}\big\la\psi(t),u\big\ra\;\textrm{ for a.e. }\;t\in[0,T]
\end{equation}
provided that the set $U$ is convex. Furthermore, for a.e.\ $t\in[0,T]$ including $t=T$ and for all $j=1,\ldots,s$ we have the
{\sc complementarity conditions}
\begin{equation}\label{41}
\big\la x^j_*,\ox(t)\big\ra<c_j\Longrightarrow\eta^j(t)=0\;\textrm{ and }\;\eta^j(t)>0\Longrightarrow\big\la x^j_*,q(t)\big\ra=c_j.
\end{equation}

\item[\bf(ii)] The {\sc transversality conditions} at the right endpoint:
\begin{equation}\label{42}
-p(T)-\sum_{j\in I(\ox(T))}\eta^j(T)x^j_*\in\lm\partial\vph\big(\ox(T)\big)\;\textrm{ and }\;\sum_{j\in I(\ox(T))}\eta^j(T)x^j_*\in N\big(\ox(T);C\big).
\end{equation}

\item[\bf(iii)] The {\sc measure nonatomicity condition:} If $t\in[0,T)$ and $\la x^j_*,\ox(t)\ra<c_j$ for all $j=1,\ldots,s$, then there
is a neighborhood $V_t$ of $t$ in $[0,T]$ such that $\gg(V)=0$ for all the Borel subsets $V$ of $V_t$.

\item[\bf(iv)] {\sc Nontriviality conditions:} It always holds that
\begin{equation}\label{e:83}
\lm+\|p(T)\|+\|q(0)\|>0.
\end{equation}
Assuming in addition that $\la x^j_*,x_0\ra<c_j$ for all $j=1,\ldots,s$, we have the {\sc enhanced nontriviality}
\begin{equation}\label{enh1}
\lm+\n p(T)\en>0.
\end{equation}
\end{itemize}
\end{theorem}\vspace*{-0.07in}
{\bf Proof.} Given the local minimizer $(\ox(\cdot),\ou(\cdot))$ for $(P)$, construct the discrete-time problems $(P_m)$ for which optimal
solutions $(\ox_m(\cdot),\ou_m(\cdot))$ exist by Proposition~\ref{Thm3.1} and converge to $(\ox(\cdot),\ou(\cdot))$ in the sense of
Theorem~\ref{ThmStrong}. We derive each of the claimed necessary conditions in $(P)$ by passing to the limit from those in
Theorem~\ref{Thm5.2*}. Let us split the derivation into several steps.\\[1ex]
{\bf Step~1:} {\em Verifying the primal equation and complementarity condition.} First we prove \eqref{37} together with the first
complementarity condition in \eqref{41}. Based on \eqref{theta}, define the functions
$$
\th_m(t):=\frac{\th^{i}_m}{h_m}\;\textrm{ for }\;t\in [t^i_m,t^{i+1}_m)\textrm{ and }\;i=0,\ldots,2^m-1
$$
on $[0,T]$ whenever $m\in\N$. It is easy to see that
\begin{eqnarray*}
\int_0^T\n\th^y_m(t)\en^2dt&=&\sum_{i=0}^{2^m-1}\frac{\Big\|\th^{iy}_m\Big\|^2}{h_m}\le\frac{1}{h_m}\sum_{i=0}^{2^m-1}
\Big(\int_{t^i_m}^{t^{i+1}_m}\n\dot{\ox}(t)-
\frac{\ox^{i+1}_m-\ox^i_m}{h_m}\en dt\Big)^2\\
&\le&\sum_{i=0}^{2^m-1}\int_{t^i_m}^{t^{i+1}_m}\Big\|\dot{\ox}(t)-\frac{\ox^{i+1}_m-\ox^i_m}{h_m}\Big\|^2dt=
\int_0^T\n\dot{\ox}(t)-\dot{\ox}_m(t)\en^2dt.
\end{eqnarray*}
Using the strong convergence $(\ox_m(\cdot),\ou_m(\cdot))\to(\ox(\cdot),\ou(\cdot))$ in Theorem~\ref{ThmStrong} ensures that
\begin{equation}\label{c:6.14}
\int_0^T\n\th^y_m(t)\en^2dt\le\int_0^T\n\dot{\ox}(t)-\dot{\ox}_m(t)\en^2dt\to0\;\textrm{ as }\;m\to\infty.
\end{equation}
This implies that a subsequence of $\{\th^y_m(t)\}$ converges, without relabeling, to zero a.e.\ on $[0,T]$. Likewise
\begin{eqnarray*}
\int_0^T\Big\|\th^u_m(t)\Big\|^2dt&=&\sum_{i=0}^{2^m-1}\frac{\Big\|\th^{iu}_m\Big\|^2}{h_m}\le\frac{1}{h_m}\sum_{i=0}^{2^m-1}
\Big(\int_{t^i_m}^{t^{i+1}_m}\n\ou^i_m-\ou(t)\en dt\Big)^2\\
&\le&\sum_{i=0}^{2^m-1}\int_{t^i_m}^{t^{i+1}_m}\n\ou^i_m-\ou(t)\en^2dt=\int_0^T\n\ou_m(t)-\ou(t)\en^2dt,
\end{eqnarray*}
which tells us, again by using Theorem~\ref{ThmStrong}, that
\begin{equation}\label{c:6.14'}
\int_0^T\n\th^u_m(t)\en^2dt\le\int_0^T\n\ou_m(t)-\ou(t)\en^2dt\to 0\;\textrm{ as }\;m\to\infty,
\end{equation}
and so $\th^u_m(t)\to 0$ for a.e.\ $t\in[0,T]$ along a subsequence. The assumed LICQ along $\ox(\cdot)$ and the robustness of this condition
yields by the choice of $z^{ji}_m$ and the convergence in Theorem~\ref{ThmStrong} that the vectors $\{z^{ji}_m\;|\;j\in I(\ox^i_m)\}$ are linearly
independent for each $i=1,\ldots,2^m$ and $m\in\N$ sufficiently large.

Taking $\eta^{i}_m\in\R^s_+$ from Theorem~\ref{Thm6.2}, we construct the piecewise constant functions $\eta_m(\cdot)$ on $[0,T]$ by
$\eta_m (t):=\eta^{i}_m $ for
$t\in [t^i_m,t^{i+1}_m)$ with $\eta_m(T):=\eta^{2^m}_m$. It follows from \eqref{87} that
\begin{equation}\label{c:51}
-\dot{\ox}_m(t)=\sum_{j=1}^s\eta^j_m(t)z^{ji}_m-g\big(\ox_m(t^i_m),\ou_m(t^i_m)\big)\;\textrm{ whenever }\;t\in (t^i_m,t^{i+1}_m),\quad m\in\N.
\end{equation}
Furthermore, we get $-\dot{\ox}(t)\in G(\ox(t))-g(\ox(t),\ou(t))$ for a.e.\ $t\in[0,T]$ with the mapping $G(\cdot)=N(\cdot;C)$, which is measurable by
\cite[Theorem~4.26]{rw}. The well-known measurable selection result (see, e.g., \cite[Corollary~4.6]{rw}) allows us to find nonnegative
measurable functions $\eta^j(\cdot)$ on $[0,T]$ for $j=1,\ldots,s$ such that equation \eqref{37} holds. Combining \eqref{c:51} and \eqref{37}
implies that
\begin{equation*}
\dot{\ox}(t)-\dot{\ox}_m(t)=\sum_{j=1}^s\big[\eta^j_m(t)z^{ji}_m-\eta^j(t)x^j_*\big]+g\big(\ox(t),\ou(t)\big)-g\big(\ox_m(t^i_m),\ou_m(t^i_m)\big)
\end{equation*}
for $t\in(t^i_m,t^{i+1}_m)$ and $i=0,\ldots,2^m-1$. It follows from the imposed LICQ
that the functions $\eta^j_m(t)$ and $\eta^j(t)$ are uniquely defined for a.e.\ $t\in[0,T]$ and belong to $L^2([0,T];\R_+)$. The constructions above yield the estimate
\begin{equation*}
\Big\|\sum_{j=1}^s\big[\eta^j(t)x^j_*-\eta^j_m(t)z^{ji}_m\big]\Big\|_{L^2}\le\n\dot{\ox}_m(t)-\dot{\ox}(t)\en_{L^2}+\n g\big(\ox(t),\ou(t)\big)-
g\big(\ox_m(t),\ou_m(t)\big)\en_{L^2}
\end{equation*}
whenever $t\in(t^i_m,t^{i+1}_m)$. Passing to the limit therein with the usage of Theorem~\ref{ThmStrong} gives us
\begin{equation*}
\sum_{j\in I(\ox(t))}\big[\eta^j(t)x^j_*-\eta^j_m(t)z^{ji}_m\big]\to 0\;\textrm{ as }\;m\to\infty\;\textrm{ for a.e. }\;t\in[0,T]
\end{equation*}
and ensures the a.e.\ convergence $\eta_m(t)\to\eta(t)$ on $[0,T]$ by the imposed LICQ. We also have that the sequence $\{\eta^{2^m}_m\}$
converges to the well-defined vector $(\eta^{1}(T),\ldots,\eta^{2^m}(T))$. Then the first complementarity condition in \eqref{41} follows from
\eqref{eta} and \eqref{eta1}.\\[1ex]
{\bf Step~2:} {\em Continuous-time extensions of approximating dual elements.} In the notation of Theorem~\ref{Thm5.2*}, define $q_m(t)$ by
extending $p^i_m(t)$ piecewise linearly on $[0,T]$ with $q_m(t^i_m):=p^i_m$ for $i=0,\ldots,2^m$.
Construct further $\gg_m(t)$ and $\psi_m(t)$ on $[0,T]$ by
\begin{equation}\label{c:6.25}
\gg_m(t):=\gg^i_m,\quad\psi_m(t):=\frac{1}{h_m}\psi^i_m\;\textrm{ for }\;t\in[t^i_m,t^{i+1}_m)\;\textrm{ and }\;i=0,\ldots,2^m-1
\end{equation}
with $\gg_m(T):=0$ and $\psi_m(T):=0$. Define now the function
\begin{equation*}
\nu_m(t):=\max\big\{t^i_m\big|\;t^i_m\le t,\;0\le i\le 2^m-1\big\}\;\textrm{ for all }\;t\in[0,T],\quad m\in\N,
\end{equation*}
and deduce respectively from $(\ref{conx})$ and $(\ref{cony})$ that
\begin{equation}\label{conx'}
\begin{array}{ll}
\dot{q}_m(t)=&-\nabla_x g\big(\ox_m(\nu_m(t)),\ou_m(\nu_m(t))\big)^*\big(-\lm_m\th^{y}_m(t)+q_m(\nu_m(t)+h_m)\big)\\\\
&+\disp\sum_{j\in I_0(-\lm_m\th^y_m(t)+q_m(\nu_m(t)+h_m ))\cup I_>(-\lm_m \th^y_m(t)+q_m(\nu_m(t)+h_m))}\gg^{j}_m(t)z^{ji}_m,
\quad\;\textrm{ and }
\end{array}
\end{equation}
\begin{equation}\label{cony1}
-\lm_m\th^{u}_m (t)-\psi_m(t)=-\nabla_u g\big(\ox_m(\nu_m(t)),\ou_m(\nu_m(t))\big)^*\big(-\lm_m\th^{y}_m(t)+q_m(\nu_m(t)+h_m)\big)
\end{equation}
for every $t\in(t^i_m,t^{i+1}_m)$ and $i=0,\ldots,2^m-1$. Next we extend the adjoint arcs $p_m(\cdot)$ to $[0,T]$ by
\begin{equation}\label{c:6.29}
p_m(t):=q_m(t)+\int_t^T\Big(\sum_{j=1}^s\gg^j_m(\tau)z^{ji}_m\Big)d\tau\;\textrm{ for every }\;t\in[0,T].
\end{equation}
This shows that $p_m(T)=q_m(T)$ and that
\begin{equation}\label{c:6.30}
\dot{p}_m(t)=\dot{q}_m(t)-\sum_{j=1}^s\gg^j_m(t)z^{ji}_m\;\textrm{ a.e. }\;t\in[0,T].
\end{equation}
The latter implies due to \eqref{conx'}, \eqref{cony1}, and the index definitions in \eqref{c56} that
\begin{equation}\label{c:59}
\dot{p}_m(t)=-\nabla_x g\big(\ox_m(\nu_m(t)),\ou_m(\nu_m(t))\big)^*\big(-\lm_m\th^{y}_m(t)+q_m(\nu_m(t)+h_m)\big)
\end{equation}
for every $t\in(t^i_m,t^{i+1}_m)$ and $i=0,\ldots,2^m-1$. Define now the vector measures $\gg^{mes}_m$ on $[0,T]$ by
\begin{equation}\label{c:6.34}
\underset{B}{\int}d\gg^{mes}_m:=\underset{B}{\int}\sum_{j=1}^s\gg^j_m(t)z^{ji}_m dt
\end{equation}
for every Borel subset $B\subset [0,T]$ and then drop for simplicity the index ``$mes$" in what follows if no confusion arises. Since all the  expressions in the statement of Theorem~\ref{Thm5.2*} are positively homogeneous of degree one with respect to $\lm_m$, $p_m$, $\gg_m$,  and $\psi_m$, the enhanced nontriviality condition \eqref{entc} and the constructions above allow us to normalize them by imposing the  sequential equality
\begin{equation}\label{c:6.35}
\lm_m +\n p_m(T)\en+\n q_m(0)\en+\int_0^T\Big\|\sum_{j=1}^s\gg^j_m(t)z^{ji}_m\Big\|dt+\int_0^T\n\psi_m(t)\en dt=1,\quad m\in\N,
\end{equation}
which tells us, in particular, that all the terms in \eqref{c:6.35} are uniformly bounded.\\[1ex]
{\bf Step~3:} {\em Verifying the dual dynamic relationships and the maximization condition.} By \eqref{c:6.35}, suppose without loss of generality
that $\lm_m\to\lm$ as $m\to\infty$ for some $\lm\ge 0$. To prove the uniform boundedness of the sequence $\{p^{0}_m,\ldots,p^{2^m}_m\}_{m\in\N}$
for all $i=0,\ldots, 2^m-1$, $m\in\mathbb{N}$, observe first from \eqref{conx} that
\begin{equation*}
p^{i+1}_m=p^{i}_m-h_m\nabla_x g(\ox^i_m,\ou^i_m)^*\Big(-\frac{1}{h_m}\lm_m\th^{iy}_m+p^{i+1}_m\Big)+h_m\sum_{j=1}^{s}\gg^{ij}_m z^{ji}_m
\end{equation*}
for all $i=0,\ldots,2^m-1$. This implies that
\begin{equation*}
\begin{aligned}
\|p^{i}_m\|&\le\|p^{i+1}_m\|+h_m\|\nabla_x g(\ox^i_m,\ou^i_m)^*\|\cdot\Big\|\Big(-\frac{1}{h_m}\lm_m\th^{iy}_m+p^{i+1}_m\Big)\Big\|+h_m
\Big\|\sum_{j=1}^{s}\gg^{ij}_m z^{ji}_m\Big\|\\
&=\big(1+h_m\|\nabla_x g(\ox^i_m,\ou^i_m)^*\|\big)\|p^{i+1}_m\|+h_m\lm_m\|\th^{y}_m(t^i_m)\|\cdot\|\nabla_x g(\ox^i_m,\ou^i_m)^*\|+h_m
\Big\|\sum_{j=1}^{s}\gg^{ij}_m z^{ji}_m\Big\|
\end{aligned}
\end{equation*}
whenever $i=0,\ldots,2^m-1$. It follows from \eqref{c:6.14} and \eqref{c:6.35} that the quantities $\nabla_x g(\ox^i_m,\ou^i_m)$,
$\lm_m\th^{iy}_m$, and $\sum_{j=1}^{s}\gg^{ij}_m z^{ji}_m$ are uniformly bounded for $i=0,\ldots,2^m-1$.
Thus we find a constant $M_1>0$ such that
\begin{equation*}
h_m\lm_m\n\th^{y}_m(t^i_m)\en\cdot\left\|\nabla_x g(\ox^i_m,\ou^i_m)^*\right\|\le M_1 h_m\n\th^{y}_m(t^i_m)\en=
M_1\sqrt{h_m\int_{t^i_m}^{t^{i+1}_m}\n\th^y_m (t)\en^2dt}
\end{equation*}
for all $i=0,\ldots,2^m-1$ and $m\in\N$. It implies that
$$
\sum_{i=0}^{2^m-1}h_m\lm_m\n\th^{y}_m(t^i_m)\en\cdot\left\|\nabla_x g(\ox^i_m,\ou^i_m)^*\right\|\le
M_1\sqrt{\int_{0}^{T}\n\th^y_m(t)\en^2dt}\to 0\;\textrm{ as }\;m\to\infty.
$$
On the other hand, we get due to \eqref{c:6.35} that
\begin{equation}\label{2ndterm}
\sum_{i=0}^{2^m-1}h_m\Big\|\sum_{j=1}^s\gg^{ij}_m z^{ji}_m\Big\|=\int_0^T\Big\|\sum_{j=1}^s\gg^j_m(t)z^{ji}_m\Big\|dt\le 1.
\end{equation}
Considering now the numbers
\begin{equation*}
A^i_m:=h_m\lm_m\n\th^{y}_m(t^i_m)\en\cdot\n\nabla_x g(\ox^i_m,\ou^i_m)^*\en+h_m\Big\|\sum_{j=1}^{s}\gg^{ij}_m z^{ji}_m\Big\|
\end{equation*}
for $i=0,\ldots,2^m-1$ and using the aforementioned uniform boundedness, find a constant $M_2>0$ such that
$\sum_{i=0}^{2^m-1}A^i_m\le M_2$. Combining the latter with the estimates above tells us that
\begin{equation}\label{t:67}
\|p^{i}_m\|\le\big(1+M_1h_m\big)\|p^{i+1}_m\|+A^i_m,\quad i=0,\ldots,2^m-1.
\end{equation}
Proceeding further by induction, we get the inequalities
\begin{eqnarray*}
\|p^{i}_m\|&\le&\big(1+M_1h_m\big)^{2^m-i}\|p^{2^m}_m\|+\sum_{j=i}^{2^m-1}A^j_m(1+M_1h_m)^{j-i}\\
&\le& e^{M_1}+e^{M_1}\sum_{i=0}^{2^m-1} A^i_m\le e^{M_1}(1+M_2)\;\textrm{ for }\;i=2,\ldots,2^m-1,
\end{eqnarray*}
which imply in turn the estimate
\begin{equation*}
\|p^{i}_m\|\le M_3\;\textrm{ for some }\;M_3>0\;\textrm{ and all }\;i=2,\ldots,2^m-1.
\end{equation*}
Hence the boundedness of $\{p^{0}_m\}$ and $\{p^{1}_m\}$ follows from \eqref{t:67} and the boundedness of $\{p^{i}_m\}_{2\le i\le 2^m}$,
which thus justifies the boundedness of the whole bundle $\{( p^{0}_m,\ldots,p^{2^m}_m)\}_{m\in\N}$.

To verify the uniform boundedness of $q_m(\cdot)$, derive from their constructions and \eqref{conx} that
\begin{equation}\label{t:68}
\begin{array}{ll}
\disp\sum_{i=0}^{2^m-1}\n q_m(t^{i+1}_m)-q_m(t^i_m)\en&\disp\le
h_m\sum_{i=0}^{2^m-1}\|\nabla_x g(\ox^i_m,\ou^i_m)^*(-\lm_m\th^{y}_m(t^i)+p^{i+1}_m)\|\\
&+\disp\int_0^T\Big\|\sum_{j=1}^s\gg^{j}_m(t)z^{ji}_m\Big\|dt
\end{array}
\end{equation}
and observe furthermore that
$$
h_m\sum_{i=0}^{2^m-1}\|\nabla_x g(\ox^i_m,\ou^i_m)^*(-\lm_m\th^{y}_m(t^i)+p^{i+1}_m)\|\le T\underset{0\le i\le
2^m-1}\max\big\{\|\nabla_x g(\ox^i_m,\ou^i_m)^*(-\lm_m\th^{y}_m(t^i)+p^{i+1}_m)\|\big\}.
$$
The latter ensures the boundedness of the first term on the right-hand side of \eqref{t:68} due to the boundedness of $\{p^{i}_m\}_{m\in\N}$,
while the boundedness of the second term therein follows from \eqref{2ndterm}. Thus we get from \eqref{t:68} that the functions $q_m(\cdot)$ on
$[0,T]$ are of uniform bounded variation on $[0,T]$ and that
\begin{equation*}
2\n q_m(t)\en-\n q_m(0)\en-\n q_m(T)\en\le\n q_m(t)-q_m(0)\en+\n q_m(T)-q_m(t)\en\le\textrm{var}(q_m;[0,T])
\end{equation*}
for all $t\in[0,T]$. Thus the sequence $\{q_m(\cdot)\}$ is bounded on $[0,T]$ since the boundedness of $\{q_m(0)\}$ and $\{q_m(T)\}$ follows
from  \eqref{c:6.35}. Applying now Helly's selection theorem gives us a function of bounded variation $q(\cdot)$ such that $q_m(t)\to q(t)$ as
$m \to\infty$ pointwise on $[0,T]$.

We see from \eqref{c:6.34} and \eqref{c:6.35} that the measure sequence $\{\gg_m\}$ is bounded in $C^*([0,T];\R^n)$. Thus the weak$^*$
sequential compactness of bounded sets in this space allows us to find a measure $\gg\in C^*([0,T];\R^n)$ such that $\{\gg_m\}$ weak*
converges to $\gg$ in $C^*([0,T];\R^n)$ along a subsequence. It follows from \eqref{c:59}, \eqref{c:6.35}, and the uniform boundedness of
$q_m(\cdot)$ on $[0,T]$ that the sequence $\{p_m(\cdot)\}$ is bounded in $W^{1,2}([0,T];\R^{n})$ and thus weakly compact in this space.
By Mazur's theorem we conclude that a sequence of convex combinations of $\dot{p}_m(\cdot)$ converges to some
$\dot{p}(\cdot)\in W^{1,2}([0,T];\R^{n})$ a.e.\ pointwise on $[0,T]$. This gives us \eqref{c:6.6} by passing to the limit along
\eqref{c:59} as $m\to\infty$ with the usage of \eqref{c:6.14} and \eqref{c:6.14'}. Note also that
\begin{equation}
\Big\|\int_t^T\sum_{j=1}^s\gg^j_m(\tau)z^{ji}_m d\tau-\int_{(t,T]}d\gg(\tau)\Big\|=\Big\|\int_t^T d\gg_m(\tau)-\int_{(t,T]}d\gg(\tau)\Big\|\to 0\;\textrm{ as }\;m\to\infty
\end{equation}
for all $t\in[0,T]$ except a countable subset of $[0,T]$ by the weak$^*$ convergence of the measures $\gg_m$ to $\gg$ in
$C^*([0,T];\R^n)$; cf.\  \cite[p.\ 325]{v} for similar arguments. Hence we get the convergence5
\begin{equation}
\int^T_t\sum_{j=1}^s\gg^j_m(\tau)z^{ji}_m d\tau\to\int_{(t,T]}d\gg(\tau)\;\textrm{ on }\;[0,T]\;\textrm{ as }\;m\to\infty
\end{equation}
and thus arrive at \eqref{c:6.9} by passing to the limit in \eqref{c:6.29}. The second (dual) complementarity condition in \eqref{41} follows from
\eqref{96} while arguing by contradiction with the usage of the established a.e.\ pointwise convergence of the functions involved therein.

To finish the proof of (i), it remains verifying the validity of the inclusion in \eqref{c:6.6'} and the maximization condition \eqref{max}. Using the
strong convergence of the discrete optimal solutions from Theorem~\ref{ThmStrong}, the convergence of $(\th^y_m(t),\th^u_m(t))\to(0,0)$ for a.e.
\ $t\in[0,T]$ obtained above as well as the robustness of the normal cone \eqref{c54}, we arrive at \eqref{c:6.6'} by passing the limit in \eqref
{cony} and in the inclusions $\psi^{i}_m\in N(\ou^i_m;U)$, $i=0,\ldots,2^m-1$, of Theorem~\ref{Thm6.2}. If $U$ is convex, the maximization
condition \eqref{max} follows directly from \eqref{c:6.6'} due to the structure \eqref{nc} of the normal cone to convex sets.\\[1ex]
{\bf Step~4:} {\em Verifying transversality inclusions.} It follows from \eqref{nmutx} and representation \eqref{F} that
\begin{equation}\label{c:72}
-p^{2^m}_m-\lm_m\vth^{2^m}_m=\sum_{j=1}^s\eta^{2^mj}_m z^{j2^{m}}_m=\sum_{j\in I(\ox^{2^m}_m)}\eta^{2^mj}_m z^{j2^{m}}_m
\in N(\ox^{2^m}_m;C^{2^m}_m),
\end{equation}
where $\eta^{2^mj}_m=0$ for $j\in\{1,\ldots,s\}\setminus I(\ox^{2^m}_m)$. Denoting
$\zeta_m:=\sum_{j\in I(\ox^{2^m}_m)}\eta^{2^mj}_m z^{j2^{m}}_m$, observe that a subsequence $\{\zeta_m\}$ converges to some
$\zeta\in\R^n$ due to the boundedness of $\lm_m$ by \eqref{c:6.35}
and the convergence of $\{p^{2^m}_m\}$ and $\{\ox^{2^m}_m \}$ with taking into account the robustness of the subdifferential. It follows from the
robustness of the normal cone in \eqref{c:72}, the convergence of $\ox^{2^m}_m\to \ox(T)$, and the inclusion $I(\ox^{2^m}_m)\subset I(\ox(T))$
for all $m$ sufficiently large that $\zeta\in N(\ox(T);C)$. Thus we get from \eqref{nmutx}) that
\begin{equation*}
-p^{2^m}_m-\zeta_m\in\lm_m\partial\vph(\ox^{2^m}_m)\;\textrm{ for all }\;m\in\N.
\end{equation*}
Passing now to the limit therein as $m\to\infty$ verifies both inclusions in \eqref{42}.\\[1ex]
{\bf Step~5:} {\em Verifying measure nonatomicity.} Take $t\in[0,T]$ with $\la x^j_*,\ox(t)\ra<c_j$ for all $j=1,\ldots,s$ and by continuity of
$\ox (\cdot)$ find a neighborhood $V_t$ of $t$ such that $\la x^j_*,\ox(\tau)\ra<c_j$ whenever $\tau\in V_t$ and $j=1,\ldots,s$.
Invoking Theorem~\ref{ThmStrong} tells us that
$\la z^{ji}_m,\ox_m(t^i_m)\ra<c^{ji}_m$ if $t^i_m\in V_t$ for all $j=1,\ldots,s$ and $m\in\N$ sufficiently large. Then we deduce from \eqref{94}
that $\gg_m(t)=0$ on any Borel subset $V$ of $V_t$. Hence
\begin{equation}\label{nonatom}
\|\gg_m\|(V)=\disp\int_Vd\|\gg_m\|=\int_V\|\gg_m(t)\|dt=0
\end{equation}
by the construction of $\gg_m$ in \eqref{c:6.34}. Passing now to limit therein and taking into account the measure convergence established
above, we get
$\|\gg\|(V)=0$, which justifies the claimed measure nonatomicity.\\[1ex]
{\bf Step~6:} {\em Verifying nontriviality conditions.} First we establish \eqref{e:83} under the general assumptions of the theorem. Arguing by
contradiction, suppose that $\lm=0$, $p(T)=0$ and $q(0)=0$. Thus $\lm_m\to 0$, $p_m(T)\to 0$, and $q_m(0)\to 0$ as $m\to\infty$. It follows
from \eqref{c:6.25} that
\begin{equation}\label{mx}
\int_0^T\n\gg_m(t)\en dt=\sum_{i=0}^{2^m-1}h_m\n\gg^i_m\en\;\textrm{ and }\;\int_0^T\n\psi_m(t)\en dt=
\sum_{i=0}^{2^m-1}h_m\frac{\n\psi^i_m \en}{h_m}=\sum_{i=0}^{2^m-1}\n\psi^i_m\en.
\end{equation}
Let us now verify the limiting condition
\begin{equation}\label{gg-lim}
\int_0^T\Big\|\sum_{j=1}^s\gg^j_m(t)z^{ji}_m\Big\|dt\to 0\;\textrm{ as }\;m\to\infty.
\end{equation}
Indeed, by $q_m(T)=p_m(T)$ and the assumption above we get $q_m(T)\to 0$ and thus deduce from \eqref{94} and \eqref{mx} that
$\int_0^T\| \gg_m(t)\|dt\to 0$ as
$m\to\infty$. Recalling that $\gg^{ij}_m=0$ for $i=0,\ldots,2^m-1$ and $j=1,\ldots,s$ by \eqref{94} and remembering the
weak$^*$ convergence of $\gg_m(\cdot)\to\gg(\cdot)$ in $C^*([0,T];\R^n)$ yield $\dot{p}(t)=\dot{q}(t)$ for a.e.\ $t\in[0,T]$ by passing to the limit in
\eqref{c:6.30}. Thus \eqref{c:6.6} reduces in this case to the linear ODE
$$
\dot{q}(t)=-\nabla_x g\big(\ox(t),\ou(t)\big)^*q(t)\;\textrm{ with }\;q(0)=0,
$$
which has only the trivial solution $q(t)\equiv 0$ on $[0,T]$. This implies that
\begin{equation}\label{max-p}
\underset{i=0,\ldots,2^m}\max\big\{\|p^{i}_m\|\big\}=\underset{t\in[0,T]}\max\big\{\|q_m(t)\|\big\}\to 0\;\textrm{ as }\;m\to\infty.
\end{equation}
By the constructions above we can estimate the left-hand side of \eqref{gg-lim} by
\begin{eqnarray*}
&&\disp\int_0^T\Big\|\sum_{j=1}^s\gg^j_m(t)z^{ji}_m\Big\|dt=\sum_{i=0}^{2^m-1}\Big\|h_m\sum_{j=1}^s\gg^{ij}_m z^{ji}_m\Big\|\\
&&\le\disp\sum_{i=0}^{2^m-1}\Big\|p^{i+1}_m-p^{i}_m+
h_m\nabla_x g(\ox^i_m,\ou^i_m)^*\Big(-\frac{1}{h_m}\lm_m\th^{iy}_m+p^{i+1}_m\Big)\Big\|\\
&&\le\disp\sum_{i=0}^{2^m-1}\big\|p^{i+1}_m\big\|\cdot\big\|\big(1+h_m\nabla_x g(\ox^i_m,\ou^i_m)^*\big)\big\|+\sum_{i=0}^{2^m-1}\n p^{i}_m
\en+\disp\sum_{i=0}^{2^m-1}\n\nabla_x g(\ox^i_m,\ou^i_m)^*\lm_m\th^{iy}_m\en.
\end{eqnarray*}
Then \eqref{max-p} and the uniform boundedness of $\nabla_x g(\ox^i_m,\ou^i_m)$ ensure that the first two terms in the last line of the obtained
estimate disappear as $m\to\infty$. To deal with the third term therein, we get by the definition of $\th^{iy}_m$ and Theorem~\ref{ThmStrong}
that
\begin{equation}\label{th-y}
\sum_{i=0}^{2^m-1}\big\|\th^{iy}_m\big\|=\sum_{i=0}^{2^m-1}\int_{t^i_m}^{t^{i+1}_m}\Big\|\frac{\ox^{i+1}_m-\ox^i_m}{h_m}-\dot{\ox}(t)\Big\|dt=
\int_0^T\big\|\dot{\ox}_m(t)-\dot{\ox}(t)\big\|dt\to 0\;\textrm{ as }\;m\to\infty,
\end{equation}
and therefore \eqref{gg-lim} is justified. To proceed further with $\psi_m^i$ in \eqref{mx}, we get by \eqref{cony} that
\begin{equation*}
\sum_{i=0}^{2^m-1}\|\psi^i_m\|\le\sum_{i=0}^{2^m-1}\|\lm_m\th^{iu}_m\|+\sum_{i=0}^{2^m-1}\big\|\lm_m\nabla_u g(\ox^i_m,\ou^i_m)^*\th^{iy}_m
\big\|+\sum_{i=0}^{2^m-1}\big\|\lm_m\nabla_u g(\ox^i_m,\ou^i_m)^*p^{i+1}_m\big\|,
\end{equation*}
which yields $\sum_{i=0}^{2^m-1}\|\psi^i_m\|\to 0$ due to \eqref{max-p}, \eqref{th-y}, and
\begin{eqnarray*}
\sum_{i=0}^{2^m-1}\|\th^{iu}_m\|=\int_0^T\big\|\ou_m(t)-\ou(t)\big\|dt\to 0\;\textrm{ as }\;m\to\infty
\end{eqnarray*}
by Theorem~~\ref{ThmStrong}. This shows that the violation of \eqref{e:83} implies the failure of \eqref{c:6.35}, a contradiction.

To complete the proof of the theorem, it remains to verify the validity of the {\em enhanced nontriviality condition} \eqref{enh1} under the
additional assumption made. Suppose on the contrary that $(\lm,p(T))=0$ while $\la x^j_*,x_0\ra<c_j$ for all $j=1,\ldots,s$. It follows from the
above arguments in the step, by using the complementarity conditions \eqref{94}, that $\dot{p}(t)=0$ for a.e.\ $t\in[0,T]$, which yields
$p(t)=p(T) =0$ on $[0,T]$. Then we get by \eqref{c:6.9} and \eqref{nonatom} that
\begin{equation}\label{qA}
q(t)=\int_{(t,T]}d\gg(\tau)=0\;\textrm{ for all }\;t\in[0,T]\setminus A,
\end{equation}
where $A\subset[0,T]$ is a countable set. Consider the two possible cases regarding \eqref{qA}:

$\bullet$ $0\notin A$, and thus $q(0)=0$.

$\bullet$ $0\in A$. In this case the measure nonatomicity condition and the fact that $A$ is at most countable
allow us to find $\tau>0$, $\tau\not\in A$, with $\int_{(0,\tau]}d\gg(t)=0$, and thus $q(0)=\int_{(\tau,T]}d\gg(t)=0$.

Hence we always have $q(0)=0$ in \eqref{qA} while showing in this way that the failure of \eqref{enh1} contradicts the validity of \eqref{e:83}
established above. $\h$\vspace*{-0.2in}

\section{Numerical Examples}\label{exa}
\setcounter{equation}{0}\vspace*{-0.1in}

In this section we consider two examples illustrating some characteristic features and strength of the necessary optimality conditions for the
sweeping control problem $(P)$ obtained in Theorem~\ref{Thm6.1*}.

Prior to dealing with specific examples, let us present the following useful assertion, which is a consequence of the measure nonatomicity
condition. \vspace*{-0.07in}

\begin{proposition}\label{claim} Assume that $\la x^*,\ox(\tau)\ra<c_j$ for all $\tau\in[t_1,t_2]$ with $t_1,t_2\in[0,T)$ and some vector
$x^*\in\{x_{\ast}^j\;|\;j=1,\ldots,s\}$, and that the measure nonatomicity condition of Theorem~{\rm\ref{Thm6.1*}} is satisfied with the measure $\gg$.
Then we have $\gamma ([t_1,t_2])=0$ and $\gamma(\{\tau\})=0$ whenever $\tau\in[t_1,t_2]$, and so $\gamma((t_1,t_2))=
\gamma([t_1,t_2))=\gamma((t_1,t_2])=0$.
\end{proposition}\vspace*{-0.07in}
{\bf Proof.} Pick any $\tau\in[t_1,t_2]$ with $\la x^*_1,\ox(t)\ra<c_j$ and find by the measure nonatomicity condition a neighborhood $V_\tau$ of
$\tau$ in $[0,T]$ such that $\gamma(V)=0$ for all the Borel subsets $V$ of $V_\tau$; in particular, $\gamma(\{\tau\})=0$.
By $[t_1,t_2]\subset \bigcup_{\tau\in[t_1,t_2]}V_\tau$ and the compactness of $[t_1,t_2]$ we find $\tau_1,\ldots,\tau_l\in[t_1,t_2]$ with
$[t_1,t_2]\subset\bigcup_{i=1}^lV_{\tau_i}$. Fix $i=1,\ldots,l-1$ and take $\Tilde{\tau_i}\in V_{\tau_i}\cap V_{\tau_{i+1}}$ with
$[\tau_i,\Tilde{\tau_i}]\subset V_{\tau_i}$ and $[\Tilde{\tau_i},\tau_{i+1}]\subset V_{\tau_{i+1}}$, where $\tau_1:=t_1$ and $\tau_l:=t_2$.
Then we arrive at the equalities
\begin{equation*}
\gamma([t_1,t_2])=\gamma\Big(\bigcup_{i=1}^{p-1}[\tau_i,\Tilde{\tau_i})\cup[\Tilde{\tau_i},\tau_{i+1})\Big)=\sum_{i=1}^{p-1}
\Big(\gamma([\tau_i,\Tilde{\tau_i}))+\gamma([\Tilde{\tau_i},\tau_{i+1}))\Big)=0,
\end{equation*}
which verify the claimed properties of the measure. $\h$\vspace*{-0.02in}

Our first example is two-dimensional with respect to both state and control variables.\vspace*{-0.07in}

\begin{example}\label{Ex-2d}
Consider the sweeping control problem of minimizing the cost functional
$$x_1(1)+x_2(1)\;\textrm{ subject to}
$$
$$\left\{\begin{matrix}
\begin{pmatrix}
\dot{x}_1\\
\dot{x}_2
\end{pmatrix}=\begin{pmatrix}
u_1\\
u_2
\end{pmatrix}-N_C\begin{pmatrix}
x_1\\
x_2
\end{pmatrix}\\
\textrm{with }\;\begin{pmatrix}
x_1\\
x_2
\end{pmatrix}(0)=\begin{pmatrix}
0\\
x_2^0
\end{pmatrix}
\end{matrix}\right.$$
where $C:=\{(x_1,x_2)\in\R^2\;|\;x_2\ge 0\}$ and $(u_1,u_2)\in U:=[-1,1]\times[-1,1]$. We rewrite the dynamics as
$$\begin{pmatrix}
\dot{x}_1\\
\dot{x}_2
\end{pmatrix}(t)=\begin{pmatrix}
u_1\\
u_2
\end{pmatrix}(t)+\eta(t)\begin{pmatrix}
0\\
1
\end{pmatrix},\quad\eta(t)\ge 0\;\textrm{ a.e. }\;t\in[0,1].$$
A direct checking shows that if $x_2^0\ge 1$ then the constraint is irrelevant and the optimal control is constant being equal to
$(-1,-1)$. If instead $0\le x_2^0<1$, then the optimal couple is $\ou_1(t)\equiv-1$ together with any measurable component $\ou_2(t)$ such that
$\ox_2(1)=0$.

The conditions of Theorem~\ref{Thm6.1*} tell us that:\\[1ex]
(1) $p=\begin{pmatrix}
p_1\\
p_2
\end{pmatrix}$ is constant on $[0,1]$ (by $(\ref{c:6.6})$);\\[1ex]
(2) $\begin{pmatrix}
-p_1\\
-p_2
\end{pmatrix}-\begin{pmatrix}
0\\
-\eta(1)
\end{pmatrix}=\begin{pmatrix}
\lambda\\
\lambda
\end{pmatrix}$,\;$\lambda\ge 0$ (by $(\ref{42})$);\\[1ex]
(3) $x^0_2>0\Longrightarrow\lambda+\|p\|>0$ (by $(\ref{enh1})$);\\[1ex]
(4) $q(t)=p-\disp\int_{[t,1]}d\gamma(\tau)=\psi(t)\in N_{[-1,1]^2}\(\begin{matrix}
\ou_1\\
\ou_2
\end{matrix}\)$ (by $(\ref{c:6.6'})$ and $(\ref{c:6.9})$);\\[1ex]
(5) $\lambda+\|p\|+\|q(0)\|>0$ (by $(\ref{e:83})$);\\[1ex]
(6)] $\eta(t)=0$ for a.e.\ $t\in[0,1]$ with $\ox_2(t)>0$\;\textrm{ and }\;\Big[$\eta(t)>0\Longrightarrow q(t)\(\begin{matrix}
0\\
1
\end{matrix}\)=0$\Big] a.e.\ $t\in[0,1]$ (by \eqref{41});\\[1ex]
(7) $d\gamma\big{|}_{\{t\:|\;\ox_2(t)>0\}}=0$ (by the measure nonatomicity condition).

To apply these conditions, consider first the case where $x_2^0>1$ in which the constraint is automatically satisfied for all the
trajectories. Since $\ox_2(1)>0$, we get $\eta(1)=0$ from (6)). If $\lambda=0$, then $p\equiv 0$ and the nontriviality condition (3) is violated. Thus we can suppose that $\lambda=1$, and so
$p=\begin{pmatrix}
-1\\
-1
\end{pmatrix}$. Condition $(7)$ implies that $d\gamma=0$ on the set in question; hence
$q\equiv p=\begin{pmatrix}
-1\\
-1
\end{pmatrix}\equiv\psi$. This shows that $\psi=\begin{pmatrix}
-1\\
-1
\end{pmatrix}$. Since $\psi\in N_{[-1,1]^2}\(\begin{matrix}
\ou_1\\
\ou_2
\end{matrix}\)$, the optimal control is $\ou(t)\equiv\begin{pmatrix}
-1\\
-1
\end{pmatrix}$.
It conforms that in this case we do not loose information with respect to the classical PMP.

Consider now the case where $0<x_2^0\le 1$. Assuming that $x_2(1)>0$ yields $\eta(1)=0$. Repeating the above arguments with the usage of (4) gives us the control
$\begin{pmatrix}
-1\\
-1
\end{pmatrix}$ on $[0,1]$ while implying that $x_2(1)=0$, a contradiction. Thus we get $x_2(1)=0$, which tells us that $\ou_2\equiv-1$ in the
case where $x_2^0=1$. Let us now deal with the first component $u_1$. Again, $\eta(1)=0$ implies that $\lambda=1$ and that $\ou_1\equiv-1$ on $[0,1]$. If $\eta(1)>0$, then $\lambda=0$ is forbidden by taking $p=-\eta(1)$ in (3), so $\ou_1\equiv-1$ is obtained as well. The case where $x_2^0=0$ requires a longer discussion, which we omit here for the sake of brevity.

This example was treated also in \cite{ac}, and the given discussion allows us to compare the two sets of necessary conditions obtained in \cite{ac} and in this paper. Actually most of them, including the adjoint equation and the transversality condition, are different. Those presented here deal only with reference trajectories where the control has bounded variation, but are more detailed and--at least in this example--are more effective for the control $u_2$ while being more difficult to use for $u_1$. This difference can be explained by the methods that are used to obtain the necessary conditions. Actually, the argument presented here takes into account the constraint at all the steps of the procedure. On the contrary, the method used in \cite{ac} is based on penalization, and so it does not see the hard constraint in the approximation steps. This explains why it behaves well with respect to $u_1$ that is not influenced by the constraint, while it is almost degenerate with respect to $u_2$.
\end{example}\vspace*{-0.05in}

The next example is also two-dimensional while addressing a more complicated polyhedral set $C$ in comparison with the halfspace in Example~\ref{Ex-2d}.\vspace*{-0.1in}

\begin{example}\label{Ex:3} Consider problem $(P)$ with the following initial data:
$$
n=m=2,\;T=1,\;x_0:=\Big(-\frac{1}{2},-\frac{1}{2}\Big),\;x^1_*:=(1,0),\;x^2_*:=(0,1),\;c_1=c_2=0,\;\varphi(x):=\frac{\n x\en^2}{2},\;g(u)=u
$$
with feasible controls $u(t)=(u^1(t),u^2(t))\in U$  a.e.\ $t\in[0,1]$ taking values in the unit square $U\subset\R^2$ with respect to the maximum
norm
$$
U:=\big\{(u^1,u^2)\in\R^2\;\big|\;\max\{u^1,u^2\}\le 1\big\}.
$$
Applying necessary optimality conditions of Theorem~\ref{Thm6.1*}, we seek for solutions to $(P)$ such that
\begin{equation}\label{ex4.3a}
\la x^j_*,\ox(t)\ra<c_j=0\;\mbox{ for all }\;t\in[0,1),\;j=1,2,\;\mbox{ and }\;\ox(1)\in{\rm bd}(C),
\end{equation}
and show that \eqref{ex4.3a} holds for $\ox(\cdot)$ found below. In the case of $(P)$ under consideration these conditions say that there exist
$\lm\ge 0$ and $\eta(\cdot)=\big(\eta^1(\cdot),\eta^2(\cdot)\big)\in L^2([0,1];\R^2_+)$ well defined at $t=1$ such that:\\[1ex]
(1) $\quad$ $\la x^*_j,\ox(t)\ra<c_j\Longrightarrow\eta^j(t)=0$ for $j=1,2$ and a.e.\ $t\in[0,1]$ including $t=1$;\\[1ex]
(2) $\quad$ $\eta^j(t)>0\Longrightarrow q^j(t)=c_j$ for $j=1,2$ and  a.e.\ $t\in[0,1]$;\\[1ex]
(3) $\quad$ $-\dot{\ox}(t)=\big(-\dot{\ox}^1(t),-\dot{\ox}^2(t)\big)=(\eta^1(t),\eta^2(t))-\big(\ou^1(t),\ou^2(t)\big)$ for a.e.\ $t\in[0,1]$;\\[1ex]
(4) $\quad$ $\big(\dot{p}^{1}(t),\dot{p}^{2}(t)\big)=\big(0,0\big)$ for a.e.\ $t\in[0,1]$;\\[1ex]
(5) $\quad$ $\big(q^{1}(t),q^{2}(t)\big)\in N\(\ou(t);U\)$ for a.e.\ $t\in[0,1]$;\\[1ex]
(6) $\quad$ $q(t)=p(t)-\gamma([t,1])$ for a.e.\ $t\in[0,1]$;\\[1ex]
(7) $\quad$ $-p(1)=\lm\big(\ox^1(1),\ox^2(1)\big)+\big(\eta^1(1),\eta^2(1)\big)$ with $\big(\eta^1(1),\eta^2(1)\big)\in N\big(\ox(1);C\big)$;\\[1ex]
(8) $\quad$ $\lm+\|p(1)\|\ne 0$.

Employing the first condition in \eqref{ex4.3a} together with (1) and (3), gives us $\dot{\ox}(t)=\ou(t)$ for a.e.\ $t\in[0,1]$. It also follows from (5) and (6) that
\begin{equation*}
q(t)=p(t)-\gg\([t,1]\)\in N\(\ou(t);U\)\;\textrm{ for a.e. }\;t\in[0,1],
\end{equation*}
which can be written in the maximization form \eqref{max}. It follows from (4) that $p(\cdot)$ is constant on $[0,1]$, i.e., $p(t)\equiv p(1)$. This
allows us to deduce that
\begin{equation*}
q(t)=p(1)-\gamma([t,1])=p(1)-\gamma(\{1\})\;\textrm{ for a.e. }\;t\in[0,1]
\end{equation*}
by using the measure nonatomicity condition of Theorem~\ref{Thm6.1*} and Proposition~\ref{claim}. Considering control functions $\ou(t)=(\vth_1,\vth_2)$ on $[0,1]$ and remembering the control constraints, we have $|\vth_1|\le 1$ and $|\vth_2|\le 1$. Thus $\ox(t)=(-\frac{1}{2}+\vth_1 t,-\frac{1}{2}+\vth_2 t)$ for all $t\in[0,1]$, and by the second condition in \eqref{ex4.3a} provides the following two possibilities:

{\bf (1)} $\ox^1(1)=0$. Then $\vth_1=\frac{1}{2}$ and the cost functional reduces is $J[\ox,\ou]=\frac{\(\vth_2-\frac{1}{2}\)^2}{2}$. It obviously
achieves its absolutely minimum value $\bar J= 0$ at the point $\vth_2=\frac{1}{2}$.

{\bf (2)} $\ox^2(1)=0$. Then $\vth_2=\frac{1}{2}$ and the minimum cost is $\bar J=0$ that is achieved at $\vth_1=\frac{1}{2}$.\\[1ex]
As a result, we arrive at a feasible solution giving the optimal value to the cost functionals:
$$
\ou(t)=\(\frac{1}{2},\frac{1}{2}\)\;\textrm{ and }\;\ox(t)=\(-\frac{1}{2}+\frac{1}{2}t,-\frac{1}{2}+\frac{1}{2}t\),\quad t\in[0,1].
$$
satisfying all the assumptions above.
\end{example}
\vspace*{-0.1in}

{\bf Acknowledgements.} The authors are grateful to Tan Cao for many useful discussions.
\end{document}